\documentclass[12pt]{article}
\usepackage{amsmath}
\usepackage{amssymb}
\usepackage{latexsym}
\newcommand{\HR}{{\cal H}_R}
\newcommand{\Ul}{{\cal U}({\cal L}^1) }
\newcommand{\M}{{\cal M}}
\newcommand{\esp}{\hspace*{0.6cm}}
\newtheorem{forme}{Theorem}[section]      
\newtheorem{corol}[forme]{Corollary}
\newtheorem{comodule}[forme]{Proposition}
\newtheorem{module}[forme]{Proposition}
\newtheorem{quotient}[forme]{Proposition}
\newtheorem{bidual}[forme]{Proposition}
\newtheorem{chaine}[forme]{Theorem}
\newtheorem{greffes}[forme]{Lemma}
\newtheorem{copperm}[forme]{Lemma}
\newtheorem{opperp}[forme]{Definition}
\newtheorem{Fi}[forme]{Lemma}
\newtheorem{decomp}[forme]{Definition}  
\newtheorem{defcom}[forme]{Proposition}
\newtheorem{th1}[forme]{Theorem}
\title{Finite dimensional comodules over the Hopf algebra of rooted trees}
\date{}
\author{L. Foissy \\
\\
{\small{\it Laboratoire de Math\'ematiques - UMR6056, Universit\'e de Reims}}\\
\small{{\it Moulin de la Housse - BP 1039 - 51687 REIMS Cedex 2, France}}\\
\small{e-mail: loic.foissy@univ-reims.fr}}

\begin{document} 
\maketitle

\section{Introduction}

 \esp In \cite{Kreimer1,Connes,Broadhurst,Kreimer2}, a Hopf algebra of rooted trees $\HR$ was introduced. 
It was shown that the antipode of this algebra was the key of a problem of renormalization (\cite{Kreimer3}).
$\HR$ is related to the Hopf algebra ${\cal H}_{CM}$ introduced in \cite{Moscovici}. Moreover, the dual algebra of $\HR$ is the enveloping algebra of the Lie algebra of rooted trees ${\cal L}^1$.
An important problem  was to give an explicit construction of the primitive elements of  $\HR$.
In \cite{Kreimer}, a  bigradation allowed to compute the dimensions of the graded parts of the space of primitive elements.

The aim of this paper is an algebraic study of $\HR$. We first use the duality theorem of \cite{Connes} to prove a result about the subcomo\-dules
of a finite dimensional comodule over the Hopf algebra of rooted trees. Then we use this result to construct comodules from 
finite families of primitive elements. Furthermore, we classify these comodules by restricting the possible families of primitive elements, and taking the quotient by the action 
of certain groups. We also show how the study of the whole algebra as a left-comodule leads to the bigrading of \cite{Kreimer}.
We then prove that ${\cal L}^1$ is a free Lie algebra.

In the next  section, we prove a formula about primitive elements of the Hopf algebra of ladders, which was already given in \cite{Kreimer}, 
and construct a projection operator on the space of primitive elements.
This operator produces the operator $S_1$ of \cite{Kreimer}.
Moreover, it allows to obtain a basis of the primitive elements by an inductive process, which answers one of the questions of \cite{Kreimer}. 

 The following sections give results about the endomorphisms of  $\HR$. First, we classify the Hopf algebra endomorphisms using the bilinear map related to the growth of trees. Then we study the coalgebra endomorphisms, 
using the graded Hopf algebra $gr(\HR)$ associated to the filtration by $deg_p$ of \cite{Kreimer}. We finally prove that $\HR \approx gr(\HR)$, and deduce a decomposition of the group of the Hopf algebra automorphisms of $\HR$ as a semi-direct product.

\section{Preliminaries}
\esp We will use notations of \cite{Connes} and \cite{Kreimer}. 
We call {\it rooted tree t} a connected and simply-connected finite set of oriented edges and vertices such that there is one distinguished vertex with no incoming edge; this vertex is called the root of $t$.
The $weight$ of $t$ is the number of its vertices. The fertility of a vertex $v$ of a tree $t$ is the number of edges outgoing from $v$. A $ladder$ is a rooted tree such that every vertex has fertility less than or equal to 1.
There is a unique ladder of weight $i$; we denote it by $l_i$.\\

\indent We define the algebra $\HR$ as the algebra of polynomials over $\mathbb{Q}$ in rooted trees. The monomials of  $\HR$ will be called $forests$.
It is often useful to think of the unit $1$ of $\HR$ as an empty forest.\\
\begin{figure}[h]
\framebox(450,60){
\begin{picture}(400,60)(-15,-10)
\circle*{5}
\put(40,0){\circle*{5}}
\put(40,0){\line(0,1){15}}
\put(40,15){\circle*{5}}
\put(80,0){\circle*{5}}
\put(80,0){\line(0,1){15}}
\put(80,15){\circle*{5}}
\put(80,15){\line(0,1){15}}
\put(80,30){\circle*{5}}
\put(130,0){\circle*{5}}
\put(130,0){\line(1,1){15}}
\put(130,0){\line(-1,1){15}}
\put(145,15){\circle*{5}}
\put(115,15){\circle*{5}}
\put(184,0){\circle*{5}}
\put(184,0){\line(0,1){15}}
\put(184,15){\circle*{5}}
\put(184,15){\line(0,1){15}}
\put(184,30){\circle*{5}}
\put(184,45){\circle*{5}}
\put(184,30){\line(0,1){15}}
\put(240,0){\circle*{5}}
\put(240,0){\line(0,1){15}}
\put(240,15){\circle*{5}}
\put(240,15){\line(1,1){15}}
\put(240,15){\line(-1,1){15}}
\put(224,30){\circle*{5}}
\put(256,30){\circle*{5}}
\put(310,0){\circle*{5}}
\put(310,0){\line(1,1){15}}
\put(310,0){\line(-1,1){15}}
\put(294,15){\circle*{5}}
\put(326,15){\circle*{5}}
\put(294,15){\line(0,1){15}}
\put(294,30){\circle*{5}}
\put(374,0){\circle*{5}}
\put(374,0){\line(1,1){15}}
\put(374,0){\line(-1,1){15}}
\put(389,15){\circle*{5}}
\put(359,15){\circle*{5}}
\put(374,0){\line(0,15){15}}
\put(374,15){\circle*{5}}
\end{picture} }
\caption{\it the rooted trees of weight less than or equal to 4. The first, second, third and fifth trees are ladders.}
\end{figure}
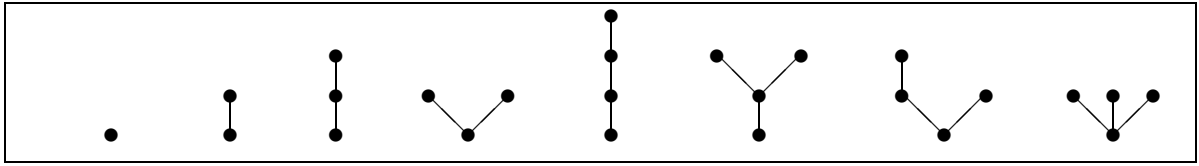

\indent We are going to give a structure of Hopf algebra to $\HR$. Before this, we define an {\it elementary cut} of a rooted tree $t$ as a cut at a single chosen edge.
An  {\it admissible cut C} of a rooted tree $t$ is an assignment of elementary cuts such that any path from any vertex of the tree has at most one elementary cut.
A cut maps a tree $t$ into a forest $t_1 \ldots t_n$. One of the $t_i$ contains the root of $t$: it will be denoted by $R^C(t)$. The product of the others will be denoted by $P^C(t)$.
Then $\Delta$ is the morphism of algebras from $\HR$ into $\HR \otimes \HR$ such that 
$$\mbox{for any rooted tree $t$, }  \Delta(t) = 1 \otimes t + t \otimes 1 + \sum_{ C \mbox{ admissible cut}}  P^C(t) \otimes R^C(t). $$

\begin{figure}[h]
{\it admissible cuts: }
\begin{picture}(27,30)(-62,0)
\put(-3,-1){\circle*{3}}
\put(-3,-1){\line(0,1){12}}
\put(-3,11){\circle*{3}}
\put(-3,11){\line(-1,1){10}}
\put(-3,11){\line(1,1){10}}
\put(-13,21){\circle*{3}}
\put(7,21){\circle*{3}}
\put(-15,12){\line(1,1){10}}
\end{picture}
\begin{picture}(27,30)(-92,0)
\put(-3,-1){\circle*{3}}
\put(-3,-1){\line(0,1){12}}
\put(-3,11){\circle*{3}}
\put(-3,11){\line(-1,1){10}}
\put(-3,11){\line(1,1){10}}
\put(-13,21){\circle*{3}}
\put(7,21){\circle*{3}}
\put(-2,21){\line(1,-1){10}}
\end{picture}
\begin{picture}(27,30)(-187,0)
\put(-3,-1){\circle*{3}}
\put(-3,-1){\line(0,1){12}}
\put(-3,11){\circle*{3}}
\put(-3,11){\line(-1,1){10}}
\put(-3,11){\line(1,1){10}}
\put(-13,21){\circle*{3}}
\put(7,21){\circle*{3}}
\put(-9,5){\line(1,0){10}}
\end{picture}
\begin{picture}(27,30)(-90,0)
\put(-3,-1){\circle*{3}}
\put(-3,-1){\line(0,1){12}}
\put(-3,11){\circle*{3}}
\put(-3,11){\line(-1,1){10}}
\put(-3,11){\line(1,1){10}}
\put(-13,21){\circle*{3}}
\put(7,21){\circle*{3}}
\put(-15,12){\line(1,1){10}}
\put(-2,21){\line(1,-1){10}}
\end{picture}

 $\Delta($
\begin{picture}(27,30)(-12,11)
\put(-3,-1){\circle*{3}}
\put(-3,-1){\line(0,1){12}}
\put(-3,11){\circle*{3}}
\put(-3,11){\line(-1,1){10}}
\put(-3,11){\line(1,1){10}}
\put(-13,21){\circle*{3}}
\put(7,21){\circle*{3}}
\end{picture}
)=  
\begin{picture}(27,30)(-12,11)
\circle*{3}
\put(-3,-1){\line(0,1){12}}
\put(-3,11){\circle*{3}}
\put(-3,11){\line(-1,1){10}}
\put(-3,11){\line(1,1){10}}
\put(-13,21){\circle*{3}}
\put(7,21){\circle*{3}}
\end{picture}
$\otimes$ 1 +
\begin{picture}(10,30)(-7,0)
\circle*{3}
\end{picture}
$\otimes$
\begin{picture}(10,40)(-2,5)
\put(-2,-1){\circle*{3}}
\put(-2,-1){\line(0,1){12}}
\put(-2,11){\circle*{3}}
\put(-2,11){\line(0,1){12}}
\put(-2,23){\circle*{3}}
\end{picture}
+
\begin{picture}(10,30)(-7,0)
\circle*{3}
\end{picture}
$\otimes$
\begin{picture}(12,40)(-2,5)
\put(-2,-1){\circle*{3}}
\put(-2,-1){\line(0,1){12}}
\put(-2,11){\circle*{3}}
\put(-2,11){\line(0,1){12}}
\put(-2,23){\circle*{3}}
\end{picture}
+ 
\begin{picture}(14,30)
\circle*{3}
\put(7,0){\circle*{3}}
\end{picture}
$\otimes$ 
\begin{picture}(10,30)(-3,2)
\circle*{3}
\put(-3,0){\line(0,1){12}}
\put(-3,12){\circle*{3}}
\end{picture}
+
\begin{picture}(27,30)(-18,11)
\put(-3,11){\circle*{3}}
\put(-3,11){\line(-1,1){10}}
\put(-3,11){\line(1,1){10}}
\put(-13,21){\circle*{3}}
\put(7,21){\circle*{3}}
\end{picture}
$\otimes$
\begin{picture}(10,30)(-2,0)
\circle*{3}
\end{picture}
+ 1 $\otimes$
\begin{picture}(27,30)(-12,11)
\put(-3,-1){\circle*{3}}
\put(-3,-1){\line(0,1){12}}
\put(-3,11){\circle*{3}}
\put(-3,11){\line(-1,1){10}}
\put(-3,11){\line(1,1){10}}
\put(-13,21){\circle*{3}}
\put(7,21){\circle*{3}}
\end{picture}

\begin{picture}(-510,-200)
\put(-5,-5){\line(0,1){100}}
\put(-5,95){\line(1,0){450}}
\put(-5,-5){\line(1,0){450}}
\put(445,-5){\line(0,1){100}}
\end{picture}
\caption{\it an example of coproduct. }
\end{figure}
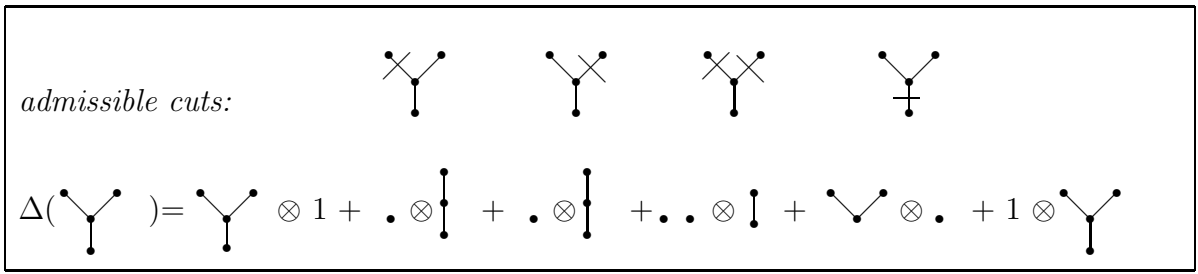 

\indent The counit is given by $\varepsilon(1)=1$, $\varepsilon(t)=0$ for any rooted tree $t$.\\
Then $\HR$ is a Hopf algebra, with antipode given by :
$$ S(t) = \sum_{\mbox{all cuts of $t$}}(-1)^{n_C+1} P^C(t) R^C(t)$$
where $n_C$ is the number of elementary cuts in $C$.\\

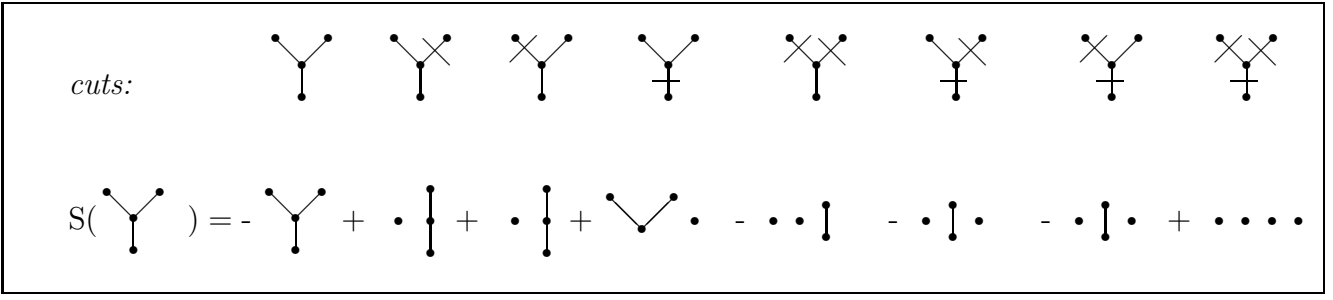
\begin{figure}[h]
{\it   cuts:}
\begin{picture}(27,30)(-62,0)
\put(-3,-1){\circle*{3}}
\put(-3,-1){\line(0,1){12}}
\put(-3,11){\circle*{3}}
\put(-3,11){\line(-1,1){10}}
\put(-3,11){\line(1,1){10}}
\put(-13,21){\circle*{3}}
\put(7,21){\circle*{3}}
\end{picture}
\begin{picture}(18,30)(-76,0)
\put(-3,-1){\circle*{3}}
\put(-3,-1){\line(0,1){12}}
\put(-3,11){\circle*{3}}
\put(-3,11){\line(-1,1){10}}
\put(-3,11){\line(1,1){10}}
\put(-13,21){\circle*{3}}
\put(7,21){\circle*{3}}
\put(-2,21){\line(1,-1){10}}
\end{picture}
\begin{picture}(27,30)(-100,0)
\put(-3,-1){\circle*{3}}
\put(-3,-1){\line(0,1){12}}
\put(-3,11){\circle*{3}}
\put(-3,11){\line(-1,1){10}}
\put(-3,11){\line(1,1){10}}
\put(-13,21){\circle*{3}}
\put(7,21){\circle*{3}}
\put(-15,12){\line(1,1){10}}
\end{picture}
\begin{picture}(27,30)(-117,0)
\put(-3,-1){\circle*{3}}
\put(-3,-1){\line(0,1){12}}
\put(-3,11){\circle*{3}}
\put(-3,11){\line(-1,1){10}}
\put(-3,11){\line(1,1){10}}
\put(-13,21){\circle*{3}}
\put(7,21){\circle*{3}}
\put(-9,5){\line(1,0){10}}
\end{picture}
\begin{picture}(27,30)(-142,0)
\put(-3,-1){\circle*{3}}
\put(-3,-1){\line(0,1){12}}
\put(-3,11){\circle*{3}}
\put(-3,11){\line(-1,1){10}}
\put(-3,11){\line(1,1){10}}
\put(-13,21){\circle*{3}}
\put(7,21){\circle*{3}}
\put(-15,12){\line(1,1){10}}
\put(-2,21){\line(1,-1){10}}
\end{picture}
\begin{picture}(27,30)(-164,0)
\put(-3,-1){\circle*{3}}
\put(-3,-1){\line(0,1){12}}
\put(-3,11){\circle*{3}}
\put(-3,11){\line(-1,1){10}}
\put(-3,11){\line(1,1){10}}
\put(-13,21){\circle*{3}}
\put(7,21){\circle*{3}}
\put(-2,21){\line(1,-1){10}}
\put(-9,5){\line(1,0){10}}
\end{picture}
\begin{picture}(27,30)(-192,0)
\put(-3,-1){\circle*{3}}
\put(-3,-1){\line(0,1){12}}
\put(-3,11){\circle*{3}}
\put(-3,11){\line(-1,1){10}}
\put(-3,11){\line(1,1){10}}
\put(-13,21){\circle*{3}}
\put(7,21){\circle*{3}}
\put(-9,5){\line(1,0){10}}
\put(-15,12){\line(1,1){10}}
\end{picture}
\begin{picture}(27,30)(-212,0)
\put(-3,-1){\circle*{3}}
\put(-3,-1){\line(0,1){12}}
\put(-3,11){\circle*{3}}
\put(-3,11){\line(-1,1){10}}
\put(-3,11){\line(1,1){10}}
\put(-13,21){\circle*{3}}
\put(7,21){\circle*{3}}
\put(-15,12){\line(1,1){10}}
\put(-2,21){\line(1,-1){10}}
\put(-9,5){\line(1,0){10}}
\end{picture}

 S(
\begin{picture}(27,30)(-13,7)
\put(-3,-1){\circle*{3}}
\put(-3,-1){\line(0,1){12}}
\put(-3,11){\circle*{3}}
\put(-3,11){\line(-1,1){10}}
\put(-3,11){\line(1,1){10}}
\put(-13,21){\circle*{3}}
\put(7,21){\circle*{3}}
\end{picture}
) = -
\begin{picture}(27,30)(-16,7)
\put(-3,-1){\circle*{3}}
\put(-3,-1){\line(0,1){12}}
\put(-3,11){\circle*{3}}
\put(-3,11){\line(-1,1){10}}
\put(-3,11){\line(1,1){10}}
\put(-13,21){\circle*{3}}
\put(7,21){\circle*{3}}
\end{picture}
 + 
\begin{picture}(27,35)(-12,-3)
\put(-3,0){\circle*{3}}
\put(9,-12){\circle*{3}}
\put(9,-12){\line(0,1){12}}
\put(9,0){\circle*{3}}
\put(9,0){\line(0,1){12}}
\put(9,12){\circle*{3}}
\end{picture}
 +
\begin{picture}(27,30)(-13,-3)
\put(-3,0){\circle*{3}}
\put(9,-12){\circle*{3}}
\put(9,-12){\line(0,1){12}}
\put(9,0){\circle*{3}}
\put(9,0){\line(0,1){12}}
\put(9,12){\circle*{3}}
\end{picture}
 +
\begin{picture}(47,30)(-10,-3)
\put(25,0){\circle*{3}}
\put(5,-3){\circle*{3}}
\put(5,-3){\line(1,1){12}}
\put(5,-3){\line(-1,1){12}}
\put(-7,9){\circle*{3}}
\put(17,9){\circle*{3}}
\end{picture}
 -
\begin{picture}(47,30)(-10,-3)
\put(-3,0){\circle*{3}}
\put(7,0){\circle*{3}}
\put(17,-6){\circle*{3}}
\put(17,6){\circle*{3}}
\put(17,-6){\line(0,1){12}}
\end{picture}
-
\begin{picture}(47,50)(-10,-3)
\put(-3,0){\circle*{3}}
\put(17,0){\circle*{3}}
\put(7,-6){\circle*{3}}
\put(7,6){\circle*{3}}
\put(7,-6){\line(0,1){12}}
\end{picture}
 -
\begin{picture}(37,50)(-10,-3)
\put(-3,0){\circle*{3}}
\put(17,0){\circle*{3}}
\put(7,-6){\circle*{3}}
\put(7,6){\circle*{3}}
\put(7,-6){\line(0,1){12}}
\end{picture}
 + 
\begin{picture}(27,50)(-10,-3)
\put(-3,0){\circle*{3}}
\put(7,0){\circle*{3}}
\put(17,0){\circle*{3}}
\put(27,0){\circle*{3}}
\end{picture}

\begin{picture}(-800,-200)
\put(-25,-10){\line(0,1){110}}
\put(-25,100){\line(1,0){500}}
\put(-25,-10){\line(1,0){500}}
\put(475,-10){\line(0,1){110}}
\end{picture}
\caption{\it the antipode.} 
\end{figure}

Let $B^+$ be the operator of $\HR$ which appends each term of a forest $t_1 \ldots t_n$ to a common root.
One can show that for every $x \in \HR$,
$$\Delta(B^+(x))= B^+(x) \otimes 1 + (Id \otimes B^+)(\Delta(x)).$$

Moreover, $\HR$ is graded as Hopf algebra by $degree(t)=weight(t)$.\\
\\
\indent For example, for all $n \in \mathbb{N}^*$, 
$$ \Delta(l_n) = 1\otimes l_n + l_n \otimes 1 + \sum_{j=1}^{n-1} l_j \otimes l_{n-j}.$$
\indent So the subalgebra of $\HR$ generated by the ladders is a Hopf subalgebra; we will denote it by ${\cal H}_{ladder}$.\\
\\
\indent We will use the Lie algebra of rooted trees ${\cal L}^1$. It is the linear span of the elements $Z_t$ indexed by rooted trees.
For $ t_1,\:t_2,\: t$ rooted trees, one defines $n(t_1,t_2;t)$ as the number of elementary cuts of $t$ such that $P^C(t)=t_1$ and $R^C(t)=t_2$.
Then the Lie bracket on ${\cal L}^1$ is given by:
$$ [ Z_{t_1},Z_{t_2}]= \sum_{t} n(t_1,t_2;t) Z_t -\sum_{t} n(t_2,t_1;t) Z_t. $$
${\cal L}^1$ is graded as Lie algebra by  $degree(Z_t)=weight(t)$. The enveloping algebra $\Ul$ is graded as Hopf algebra with the corresponding gradation (see \cite{Connes, Panaite}).

\section{Duality between $ {\cal H}_R $ -comodules and $ {\cal U}({\cal L}^1) $-modules} 

\esp We shall use the following result of \cite{Connes}:

\begin{forme}
There is a bilinear form on  $ {\cal U}({\cal L}^1) \times {\cal H}_R $ defined by: 
\begin{eqnarray*}
    <1,P(t_i)>&=&\varepsilon (P(t_i)),\\
< Z_t, P(t_i)>&=&(\frac{\partial}{\partial t}P)(0),\\
 \mbox{ and }<Z_1 Z_2,P>&=&<Z_1 \otimes Z_2 ,  \Delta (P)>. 
\end{eqnarray*}
 \end{forme}

An easy induction on $weight \:(P(t_i))$ proves the following property:
\newtheorem{grforme}[forme]{Lemma}
\begin{grforme}
\label{grform}
If $l \in {\cal U}({\cal L}^1)$ and  $P(t_i) \in \HR$ are homogeneous of different degrees,\\
then $< l,P(t_i)> =0 $.
\end{grforme}

Let $ {\cal I}_n $ be the ideal of $ {\cal H}_R $ generated by the homogeneous elements of weight greater than or equal to $n$ 
and $ {\cal J}_n $ the ideal of $ {\cal U}({\cal L}^1) $ generated by the homogeneous elements of weight greater than or equal to $n$.
Let  $ {\HR}^{*g}= \{  f \in  {\cal H}_R^* / \exists n \in \mathbb{N}, f ({\cal I}_n) = (0) \} $
and  $ {\cal U}({\cal L}^1)^{*g} = \{  f \in {\cal U}({\cal L}^1)^* / \exists n \in \mathbb{N}, f ({\cal J}_n) = (0) \} $.
One defines an algebra structure on $ {\HR}^{*g}$ by dualising the coproduct on   ${\cal H}_R$ and a coalgebra structure on 
$ {\cal U}({\cal L}^1)^{*g} $ by dualising the product of $ {\cal U}({\cal L}^1) $. Then we have the following result:

\begin{corol}
$$ \mbox{Let }\Phi:
   \left\{
   \begin{array}{ccc}
   $$ {\cal H}_R $$ & $$ \longmapsto $$ & $$ {\cal U}({\cal L}^1)^{*g}$$\\
   $$ P(t_i) $$ & $$ \longmapsto $$ & $$ \langle .,P(t_i)\rangle $$
  \end{array}
  \right. 
  \mbox{ and let }\Psi:
   \left\{
   \begin{array}{ccc}
   $$ {\cal U}({\cal L}^1) $$ & $$ \longmapsto $$ & $$ {\cal H}_R^{*g}$$\\
   $$ l $$ & $$ \longmapsto $$ & $$ \langle l,.\rangle. $$
  \end{array}
  \right.  
  $$

  Then $ \Phi $ is a coalgebra isomorphism and $ \Psi $ is  an algebra isomorphism.
\end{corol}
  
One can now dualise $ {\cal H}_R $-comodules and $ {\cal U}({\cal L}^1) $-modules.
First, we have:

\begin{comodule}
\label{com}
Let C be a $ {\cal H}_R $-comodule and $ \Delta_C $ its structure map: $C \longmapsto \HR \otimes C$. 
\\
Then $ C^* $ is a $ {\cal U}({\cal L}^1) $-module with:
$$ \forall l \in  {\cal U}({\cal L}^1), \forall f \in C^*,\forall x \in C,\:l.f(x) = \sum_{(x)} \langle l,x^{(1)} \rangle f(x^{(2)}) $$

$$  \mbox{ where } \Delta_C (x) = \sum_{(x)} x^{(1)} \otimes x^{(2)}.$$
\end{comodule}

\noindent
\noindent  {\it Proof:} classical; see \cite{Sweedler}.

\begin{module}
Let M be a $ {\cal U}({\cal L}^1) $-module.
Let  $ M^{*g} = \{ f \in M^* /\exists n \in \mathbb{N}, f ({\cal J}_n M) = (0) \} $.
Then $ M^{*g} $ is a ${\cal H}_R$-comodule with $ \Delta_M :M^{*g} \longmapsto {\cal H}_R \otimes M^{*g} $ defined by:\\
$$\forall f \in M^{*g}, \forall l \in  {\cal U}({\cal L}^1), \forall x \in M, 
\mbox{ with }\Delta_M (f)= \sum_{(f)} f^{(1)} \otimes  f^{(2)}:$$ 
$$\Delta_M (f).( l \otimes m ) = \sum_{(f)} \langle l, f^{(1)} \rangle f^{(2)}(m) = f(l.m). $$

\end{module}

\noindent  {\it Proof:} 
$$ \mbox{Let } \alpha: 
   \left\{
   \begin{array}{ccc}
   $$ {\cal U}({\cal L}^1)^{*g} \otimes M^{*g} $$ & $$ \longmapsto $$ & $$ ({\cal U}({\cal L}^1) \otimes M)^*$$\\
   $$ f \otimes g $$ & $$ \longmapsto $$ & $$ \left\{  
                                               \begin{array}{clr}
                                                $$ { \cal U}({\cal L}^1) \otimes M $$ & $$ \longmapsto $$  $$ \hspace{1cm} \mathbb{Q} $$\\
                                                $$ l \otimes m $$ & \longmapsto & \hspace{-1,3cm}$$f(l)g(m); $$
                                               \end{array}
                                               \right.
                                               $$

  \end{array}
  \right. $$
  $\alpha $ is injective. If $\mu$ is the structure map of M and $\mu^*$ its transpose ($\mu:{\cal U}({\cal L}^1)\otimes M \longmapsto M$), 
  we have to show that $Im \mu^* \subset Im\alpha $. 
  With the definition of $ M^{*g} $, one easily has: \\
  $ Im\alpha = \{ f \in ({\cal U}({\cal L}^1) \otimes M)^* / \exists n \in \mathbb{N}, 
  f ({\cal J}_n \otimes M) = (0),f(A \otimes {\cal J}_n M) = (0) \} $.\\
  Let $ f \in M^{*g}, l \otimes m \in {\cal U}({\cal L}^1) \otimes M $. $ \mu^*(f)(l \otimes m ) = f(l.m).$   
  As $ f \in M^{*g}$, clearly $ \mu^*(f)$ is in $ Im \alpha $.

\begin{quotient}
\label{quot}
Let $M_1,M_2$ two ${\cal U}({\cal L}^1)$-modules, with $M_1 \subset M_2$; there exists an injection of comodules: 
$$ (M_2/M_1)^{*g} \longmapsto M_2^{*g}.$$
\end{quotient}

\noindent  {\it Proof:} let $ p:M_2 \longmapsto M_2/M_1 $ the canonical surjection; then it is easy to see that its transpose is an injective
morphism of comodules from  $ (M_2/M_1)^{*g} $ to $M_2^{*g}$.

\begin{bidual}
\label{bid}
Let C a finite-dimensional ${\cal H}_R$-comodule. Then $C^*$ is a ${\cal U}({\cal L}^1)$-module, and
${(C^*)}^{*g} $ is the whole  ${(C^*)}^* $. Moreover C and ${(C^*)}^* $ are isomorphic ${\cal H}_R$-comodules. 
\end{bidual}
\noindent  {\it Proof:} let $ l \in {\cal U}({\cal L}^1), f \in C^*, x \in C$. 
Then $ (l.f)(x) = \sum_{(x)} \langle l,x^{(1)}\rangle,f(x^{(2)})$.\\ 
Let $k_x = max_{(x)}\left(weight(x^{(1)}) \right)+1$. If $l$ is homogeneous of weight greater than $k_x$, 
then $(l.f)(x) = 0$ (lemma \ref{grform}). As C is finite-dimensional, there exists $k \in \mathbb{N}, k\ge k_x \, \forall x \in C$,  
hence ${\cal J}_k.C^* = (0)$, and hence ${(C^*)}^{*g} ={(C^*)}^* $.
It is then easy to show that the canonical isomorphism 
between C and ${(C^*)}^*$ is a comodule isomorphism.\\
\\
We are now ready to prove the:

\begin{chaine}
\label{chain}
Let C be a finite-dimensional ${\cal H}_R$-comodule and n its dimension; then $C$ has a complete flag of comodules, that is to say:\\
$ \forall i \in \{1 \ldots n\},\exists \, C^{\,(i)} \mbox{ a subcomodule of C of dimension i, with }
C^{\,(1)} \subset \ldots \subset C^{\,(n)} = C$.
\end{chaine}

\noindent  {\it Proof:}
it is enough to exhibit a subcomodule of dimension  $n-1$.
By proposition \ref{com}, $C^*$ is a ${\cal U}({\cal L}^1)$-module, and there exists $k \in \mathbb{N}, \: {\cal J}_k.C^* = (0)$. 
Hence as a ${\cal L}^1$-module, $l.C^*=(0)$ for every $l$  in ${\cal L}^1$, homogeneous of weight greater than $n$.
So $C^*$ is in fact a module over the quotient of  ${\cal L}^1$ by the Lie ideal generated by these $l$,  
and it is clear that this quotient is a finite-dimensional nilpotent Lie algebra.
Moreover, every $l \in {\cal L}^1$ is a nilpotent endomorphism of $C^*$. 
By Engel's theorem, $C^*$ has a submodule C' of dimension 1. ${\cal J}_k.(C^*/C') = (0)$ because ${\cal J}_k.C^* =(0)$,  
so ${(C^*/C')}^{*g}= {(C^*/C')}^*$, and the dimension of this comodule is $n-1$.
By proposition \ref{bid}, C is isomorphic to ${(C^*)}^*$ which has a subcomodule of dimension $n-1$ by proposition \ref{quot}.  
\\

{\it Remark:} one can use the fact that ${\cal L}^1$ acts by zero on $C'$ (which is given by Engel's theorem), to show that the quotients $\frac{C^{(i+1)}}{C^{(i)}}$ are trivial comodules,
that is to say $\Delta(\overline{x}) = 1 \otimes \overline{x}$ $\forall \overline{x} \in \frac{C^{(i+1)}}{C^{(i)}}$.  
\section{Natural growth}
\esp Let $M,N$ be two forests of ${\cal H}_R$. We define: 
$$ M \top N=
\left\{
\begin{array}{cc}
\frac{1}{weight (N)} \sum \mbox{forests obtained by appending M to every node of N } & \mbox{if } N \neq 1 \\
 0  & \mbox{if } N=1.
\end{array}
\right.
$$
\indent We extend $.\top .$ to a bilinear map from ${\cal H}_R \times {\cal H}_R$ into ${\cal H}_R$.\\

\begin{figure}[h]
\begin{picture}(15,30)(-10,0)
\put(0,0){\circle*{3}}
\put(0,0){\line(0,1){10}}
\put(0,10){\circle*{3}}
\end{picture}
$\top$
\begin{picture}(30,30)(-10,0)
\put(0,0){\circle*{3}}
\put(0,0){\line(1,1){10}}
\put(0,0){\line(-1,1){10}}
\put(-10,10){\circle*{3}}
\put(10,10){\circle*{3}}
\end{picture}
= $\frac{1}{3} \left(
\begin{picture}(30,30)(-13,15)
\put(0,0){\circle*{3}}
\put(0,0){\line(1,1){10}}
\put(0,0){\line(-1,1){10}}
\put(-10,10){\circle*{3}}
\put(10,10){\circle*{3}}
\put(-10,10){\line(0,1){10}}
\put(-10,20){\line(0,1){10}}
\put(-10,20){\circle*{3}}
\put(-10,30){\circle*{3}}
\end{picture}
+
\begin{picture}(30,30)(-13,15)
\put(0,0){\circle*{3}}
\put(0,0){\line(1,1){10}}
\put(0,0){\line(-1,1){10}}
\put(-10,10){\circle*{3}}
\put(10,10){\circle*{3}}
\put(10,10){\line(0,1){10}}
\put(10,20){\line(0,1){10}}
\put(10,20){\circle*{3}}
\put(10,30){\circle*{3}}
\end{picture}
+
\begin{picture}(30,30)(-13,15)
\put(0,0){\circle*{3}}
\put(0,0){\line(1,1){10}}
\put(0,0){\line(-1,1){10}}
\put(-10,10){\circle*{3}}
\put(10,10){\circle*{3}}
\put(0,0){\line(0,1){10}}
\put(0,10){\line(0,1){10}}
\put(0,10){\circle*{3}}
\put(0,20){\circle*{3}}
\end{picture}
\right)$
; 
\begin{picture}(20,10)
\end{picture}
\begin{picture}(20,30)(0,0)
\put(0,0){\circle*{3}}
\put(10,0){\circle*{3}}
\end{picture}
$\top$
\begin{picture}(10,30)(-10,10)
\put(0,0){\circle*{3}}
\put(0,0){\line(0,1){10}}
\put(0,10){\circle*{3}}
\put(0,10){\line(0,1){10}}
\put(0,20){\circle*{3}}
\end{picture}
$=\frac{1}{3}
\left(
\begin{picture}(25,30)(-10,10)
\put(0,0){\circle*{3}}
\put(0,0){\line(0,1){10}}
\put(0,10){\circle*{3}}
\put(0,10){\line(0,1){10}}
\put(0,20){\circle*{3}}
\put(10,10){\circle*{3}}
\put(-10,10){\circle*{3}}
\put(0,0){\line(1,1){10}}
\put(0,0){\line(-1,1){10}}
\end{picture}
+\begin{picture}(20,30)(-10,10)
\put(0,0){\circle*{3}}
\put(0,0){\line(0,1){10}}
\put(0,10){\circle*{3}}
\put(0,10){\line(0,1){10}}
\put(0,20){\circle*{3}}
\put(10,20){\circle*{3}}
\put(-10,20){\circle*{3}}
\put(0,10){\line(1,1){10}}
\put(0,10){\line(-1,1){10}}
\end{picture}
+\begin{picture}(20,40)(-10,10)
\put(0,0){\circle*{3}}
\put(0,0){\line(0,1){10}}
\put(0,10){\circle*{3}}
\put(0,10){\line(0,1){10}}
\put(0,20){\circle*{3}}
\put(10,30){\circle*{3}}
\put(-10,30){\circle*{3}}
\put(0,20){\line(1,1){10}}
\put(0,20){\line(-1,1){10}}
\end{picture}
\right)$
;

\begin{picture}(130,0)
\end{picture}
\begin{picture}(10,30)(0,0)
\put(0,0){\circle*{3}}
\put(0,0){\line(0,1){10}}
\put(0,10){\circle*{3}}
\end{picture}
$\top$
\begin{picture}(20,30)(-5,0)
\put(0,0){\circle*{3}}
\put(10,0){\circle*{3}}
\end{picture}
$= \frac{1}{2} \left(
\begin{picture}(20,30)(-5,10)
\put(0,0){\circle*{3}}
\put(10,0){\circle*{3}}
\put(0,10){\circle*{3}}
\put(0,10){\line(0,1){10}}
\put(0,20){\circle*{3}}
\put(0,0){\line(0,1){10}}
\end{picture}
+
\begin{picture}(20,30)(-5,10)
\put(0,0){\circle*{3}}
\put(10,0){\circle*{3}}
\put(10,10){\circle*{3}}
\put(10,10){\line(0,1){10}}
\put(10,20){\circle*{3}}
\put(10,0){\line(0,1){10}}
\end{picture}
\right)$

\begin{picture}(300,-200)(0,5)
\put(0,0){\line(0,1){135}}
\put(0,0){\line(1,0){500}}
\put(500,0){\line(0,1){135}}
\put(0,135){\line(1,0){500}}
\end{picture}
\caption{\it the bilinear map $\top$.}

\end{figure}

In the following we use the notation $\tilde{\Delta}(x) = \Delta(x) - 1 \otimes x -x\otimes 1$ for every $x \in {\cal H}_R$.
We have $ Prim({\cal H}_R) = Ker(\tilde{\Delta})$.

\begin{greffes}
\label{greffe}

Let $x \in {\cal H}_R$ and $y$ be a primitive element of ${\cal H}_R$. Then we have:
$$ \tilde{\Delta}(x\top y)= x \otimes y + \sum_{(x)} x^{(1)} \otimes ( x^{(2)} \top y )$$
where $ \tilde{\Delta}(x) =\sum_{(x)} x^{(1)} \otimes  x^{(2)} $. 
\end{greffes}
\noindent  {\it Proof:} see \cite{Kreimer}, section 5.4.
\begin{opperp}
\label{defFi}
\textnormal{Let $ i \in \mathbb{N}^*$ and $ p_1, \ldots, p_i $ be primitive elements of ${\cal H}_R $.
By induction on $i$ we define $ p_i \top \ldots \top  p_1 $ by $ ( p_i \top \ldots \top p_2) \top p_1$. 
And we define:} 
$$ F_i: \left \{ 
\begin{array}{ccc} 
$$ Prim({\cal H}_R)^{\otimes i}$$ & $$\longmapsto$$ & $${\cal H}_R$$\\
 $$p_i \otimes \ldots \otimes p_1$$ & $$ \longmapsto$$ & $$p_i \top \ldots \top p_1.$$
 \end{array}
 \right.
 $$
\end{opperp}
\begin{copperm}
\label{cop}
Let $ p_1, \ldots, p_i $ be primitive elements of ${\cal H}_R $.
$$ \tilde{\Delta}( p_i \top \ldots \top p_1) =\sum_{j=1}^{j=i-1} (p_i \top \ldots \top p_{j+1})
           \otimes ( p_j \top \ldots \top p_1). $$
\end{copperm}
\noindent  {\it Proof:} by induction, using \ref{greffe}.
\\
\\
One remarks easily that $\tilde{\Delta}$ is still coassociative.
 We define $\tilde{\Delta}^0=Id_{\HR} - \eta \circ \varepsilon$, $\tilde{\Delta}^1=\tilde{\Delta}$, 
and by induction $\tilde{\Delta}^k = (\tilde{\Delta}^{k-1} \otimes Id) \circ \tilde{\Delta}$.

\begin{Fi}
\label{Fis}
Let   $ i \in \mathbb{N}^*$; then $\tilde{\Delta}^{i-1} \circ F_i = Id_{{[Prim({\cal H}_R)]}^{\otimes i}}$;  
if $ k>i-1$, $ \tilde{\Delta}^k \circ F_i =0$. Moreover, $F_i$ is injective, and the sum $(1)+\sum_{i=1}^{\infty} Im(F_i)$ is direct.
\end{Fi}
\noindent  {\it Proof:} one shows the first point by induction, using \ref{cop}. The second point is an immediate corollary. 
For the last point, let $x_0 \in \mathbb{Q},$ $x_i \in Im(F_i) \: \forall i\in \{1 \ldots n\}$, with $x_01 +x_1+\ldots +x_n =0$. 
Then $ \varepsilon(0)=x_0=0$. Moreover, $\tilde{\Delta}^{n-1}(x_1+\ldots +x_n)= \tilde{\Delta}^{n-1}(x_n)=0$. 
As $x_n=F_n(y_n)$ for a certain $y_n$, we have $y_n=0$, so $x_n=0$. One concludes by an induction on $n$.

\section{ Construction and parametrization of finite-dimensional ${\cal H}_R$-comodules}

\begin{decomp}
\textnormal{
Let $(i,j) \in {(\mathbb{N}^*)}^2,i \leq j$.
We denote $I_{i,j} :=\{i \ldots j \}$. 
A decomposition of $I_{i,j}$ is a partition of $I_{i,j} $ in connected parts.
We denote a decomposition in the following way: 
    $$I_{i_1,j_1} \ldots   I_{i_k,j_k} \mbox{ with } 
    i=i_1 \leq j_1 <i_2 \leq \ldots < i_k \leq j_k = j.$$
And we denote by ${\cal D}_{i,j}$ the set of all decompositions of  $I_{i,j}$.\\
There are $2^{j-i}$ decompositions of $I_{i,j}$.  }
\end{decomp}

\begin{defcom}
\label{defcomo}
Let $n \geq 1, (p_{i,j})_{1 \leq i \leq j \leq n }$ any family of $\frac{n(n+1)}{2}$ primitive elements of ${\cal H}_R$.
Let $C$ be a vector space of dimension $n+1$, with basis $(e_0, \ldots,e_n)$.
We define: 
\begin{eqnarray*}
 \Delta_C(e_0)&=&  1 \otimes e_0;\\
\Delta_C(e_i)&=&  \left [\sum_{j=0}^{j=i-1} \left( \sum_{ I_{i_1,j_1} \ldots   I_{i_k,j_k} \in {\cal D}_{j+1,i}}
      p_{i_k,j_k} \top \ldots \top p_{i_1,j_1} \right) \otimes e_j \right]+1 \otimes e_i.
\end{eqnarray*}
Then $(C,\Delta_C)$ is a (left) ${\cal H}_R$-comodule. We denote this comodule  by $C_{(p_{i,j})}$. 
\end{defcom}
\noindent  {\it Proof:} the axiom of counity is trivial.\\
Coassociativity: we have to show that $((\Delta \otimes Id) \circ \Delta_C)(e_i)=((Id \otimes \Delta_C) \circ \Delta_C)(e_i) \, \forall i$. 
It is trivial for $i=0$. For $i \geq 1$, we have:
\begin{eqnarray*}
((Id \otimes \Delta_C) \circ \Delta_C)(e_i)  &=& 
\sum_{j=0}^{j=i} \sum_{l=0}^{l=j} \left( \sum_{{\cal D}_{j+1,i}} p_{i_k,j_k} \top \ldots \top p_{i_1,j_1} \right) \otimes \left( \sum_{{\cal D}_{l+1,j}} 
   p_{i'_r,j'_r} \top \ldots \top p_{i'_1,j'_1} \right) \otimes e_l\\ 
& =& \sum_{l=0}^{i} \sum_{{\cal D}_{l+1,i}} \Delta(p_{i''_s,j''_s} \top \ldots \top p_{i''_1,j''_1}) \otimes e_l
     \mbox{  (by \ref{cop})}\\
&=& ((\Delta \otimes Id) \circ \Delta_C)(e_i).
\end{eqnarray*}

The following theorem gives a parametrization of the finite dimensional $\HR$-comodules by certain finite families of primitive elements:

\begin{th1}
\label{theo1}
Let $(C,\Delta_C)$ be a finite-dimensional comodule. If the dimension of $C$ is 1, then $C$ is trivial, that is to say $\Delta_C(x)=1\otimes x \: \: \forall x \in C$. 
If the dimension of C is $n,\:n \geq 2$, then there is a finite family $(p_{i,j})_{1 \leq i \leq j \leq n }$ of $\frac{n(n+1)}{2}$ primitive elements of $\HR$ such that
$C$ is isomorphic to $C_{(p_{i,j})}$. 
\end{th1}
We shall use the following lemma:
\newtheorem{groupo}[forme]{Lemma}
\begin{groupo}
If $x \in {\cal H}_R$ is such that $\Delta(x) = x \otimes x$, then $ x=0 \mbox{ or } 1$. 
\end{groupo}
\noindent  {\it Proof:} suppose $ x \neq 0$. As $\Delta$ is homogeneous of degree 0, $x$ is of weight 0. 
It is then trivial that $x=1$. \\
\\
{\it Proof of the theorem:} let $C^{\,(0)} \subset \ldots \subset C^{\,(n)} $ be a complete flag of subcomodules, which exists by \ref{chain}, and let $(e_0, \ldots ,e_n)$ be an adapted basis to this flag. 
Then we have a fa\-mily $(Q_{i,j})_{1 \leq j \leq i \leq n}$ of elements of ${\cal H}_R$ such that 
$\Delta(e_i)=\sum_{j=0}^{j=i} Q_{i,j} \otimes e_j$. (If $n=0$, then $(Q_{i,j})_{1 \leq j \leq i \leq n}$  is empty).
The axiom of counity implies that $\varepsilon(Q_{i,i})=1$, 
and $\Delta(Q_{i,j})=\sum_{l=j}^{l=i} Q_{i,l} \otimes Q_{l,j}$ by the axiom of coassociativity.
So by the lemma, $Q_{i,i}=1 \: \forall i$, which proves the theorem for $n=0$.
Moreover, $Q_{i,i-1}$ is primitive. If $n=1$, $C \approx C_{(p_{1,1})}$ with $ p_{1,1}=Q_{1,0}$. 
We end with an induction on $n$: by induction hypothesis on $C'$ spanned by $(e_0, \ldots, e_{n-1})$, 
we have $p_{i,j},\:1\leq i \leq j \leq n-1$. With $p_{n,n}=Q_{n,n-1}$, we have
$Q_{n,n-1}=\sum_{{\cal D}_{n,n}} p_{i_k,j_k} \top \ldots \top p_{i_1,j_1}$.
Suppose we have built $p_{n,n}, \ldots,p_{i+1,n}$, such that 
 $Q_{n,i}=\sum_{{\cal D}_{i+1,n}} p_{i_k,j_k} \top \ldots \top p_{i_1,j_1}$.
 Then
\begin{eqnarray*}
  \tilde{\Delta}(Q_{n,i-1})&=&\sum_{l=i}^{l=n-1} \left( \sum_{{\cal D}_{l+1,n}}
p_{i_k,j_k} \top \ldots \top p_{i_1,j_1} \right) \otimes \left( \sum_{{\cal D}_{i,l}} 
   p_{i'_r,j'r} \top \ldots \top p_{i'_1,j'_1} \right) \\ 
&=&\sum_{{\cal D}_{i,n}-\{I_{i,n}\}}  \tilde{\Delta}(p_{i''_s,j''_s} \top \ldots \top p_{i''_1,j''_1}).
\end{eqnarray*}
$$\hspace{-1.8cm}\mbox{As $Ker(\tilde{\Delta})=Prim({\cal H}_R)$, we take }
p_{i,n}=Q_{n,i-1}-\sum_{{\cal D}_{i,n}-\{I_{i,n}\}}  (p_{i''_s,j''_s}\top \ldots \top p_{i''_1,j''_1}).
$$

\noindent \textbf{5.6} {\it Remarks:} \begin{enumerate}
\item The family $(p_{i,j})$ depends on the choice of the basis $(e_0,\ldots,e_n)$, hence is not unique.
\item By the following, we shall identify $(p_{i,j})_{1 \leq i \leq j \leq n }$ with
$$\left [ \begin{array}{cccc}
    0 & p_{1,1} & \cdots & p_{1,n}\\
    \vdots & \ddots & \ddots & \vdots\\
    0 & \cdots & \ddots & p_{n,n}\\
    0 & \cdots & \cdots & 0
    \end{array}
    \right]
    = {\cal P} \in {\cal M}_{n+1}(Prim({\cal H}_R))
    $$
where  ${\cal M}_{n+1}(Prim({\cal H}_R))$ is the space of square matrices of order $n+1$ with entries in $Prim(\HR)$.    
With the notation of the proof of \ref{theo1}, we will write $$ 
{\cal Q} = 
\left [ \begin{array}{cccc}
    Q_{0,0} & 0 &\cdots & 0\\
    \vdots & \ddots & \ddots & \vdots\\
    Q_{n-1,0} & \cdots & \ddots & 0\\
    Q_{n,0} & \cdots & \cdots & Q_{n,n}
    \end{array}
    \right] \in {\cal M}_{n+1}({\cal H}_R)$$
where  ${\cal M}_{n+1}({\cal H}_R)$ is the space of square matrices of order $n+1$ with entries in $\HR$.    
\\
Recall that $F_i$ was defined in \ref{defFi}.
Let $\pi_1$ be the projection on $Prim({\cal H}_R)=Im(F_1)$ 
 in $(1) \oplus \oplus_{i=1}^{i=\infty}Im(F_i)$. Then $Q_{i,j} \in (1) \oplus \oplus_{i=1}^{i=\infty}Im(F_i)$, and 
 $\pi_1(Q_{i,j})=p_{j+1,i}$, or in a matricial form: ${\cal P}=\pi_1({\cal Q}^T)$ (here $\pi_1$ acts on each entry of the matrix).
\end{enumerate}

\section{Classification of the finite-dimensional ${\cal H}_R$-comodules}
\newtheorem{reduite}[forme]{Definition}
\begin{reduite}
\textnormal{ Let ${(p_{i,j})}_{1 \leq i \leq j \leq n}$ be a family of $\frac{n(n+1)}{2}$ primitive elements of $\HR$ and ${\cal P}$ \\[1mm]
the associated matrix as in the remark 5.6. We say that $(p_{i,j})$ is reduced
if there are $c_0,\ldots,c_k \in \mathbb{N}^*$ such that:
$${\cal P}=
\left [ \begin{tabular}{c|c|c|c}
    $0$ & ${\cal P}_{1,1}$ & $\cdots$ & ${\cal P}_{1,k}$\\
    \hline
    $\vdots$ & $\ddots$ & $\ddots$ & $\vdots$\\
    \hline
    $0$ & $\cdots$ & $\ddots$ & ${\cal P}_{k,k}$\\
    \hline
    $0$ & $\cdots$ & $\cdots$ & $0$
    \end{tabular}
    \right]
    $$
where the diagonal zero blocs are in ${\cal M}_{c_0}({\cal H}_R), \ldots ,{\cal M}_{c_k}({\cal H}_R)$ 
and the columns in each bloc ${\cal P}_{i,i}$, $1 \leq i \leq k$, are linearly independent;
$(c_0,\ldots,c_k)$ is called the type of $(p_{i,j})$.  }
\end{reduite}

{\it Example:} 
$$ \mbox{Let }{\cal P}= \left[
\begin{tabular}{c|cc|cc}
$0$ & $a$ & $b$ & $x$ & $y$\\
 \hline
0 & 0 & 0 & $c$ & $e$\\
0 & 0 & 0 & $d$ & $f$ \\
 \hline
0 & 0 & 0 & 0 & 0\\
0 & 0 & 0 & 0 & 0
\end{tabular} \right]
\in {\cal M}_5(Prim({\cal H}_R)).$$
Suppose that $a$ and $b$ are linearly independent in the vector space $\HR$, and 
$\left( \begin{array}{c} c \\ d \end{array} \right)$ and $\left( \begin{array}{c} e \\ f \end{array} \right)$ are linearly independent in the vector space $\HR^2$.
Then $(p_{i,j})$ is a reduced family of type (1,2,2).

\newtheorem{drapeau}[forme]{Definition}
\begin{drapeau}
\label{drap}
\textnormal{ Let $C$ be a ${\cal H}_R$-comodule. One defines $C_0=\{ x \in C/ \Delta_C(x) = 1 \otimes x\}$ 
and, by induction, $C_{i+1}$  the unique subcomodule of $C$ such that 
\begin{description}
\item[\it i)] $C_i \subset C_{i+1}; $\\
\item[\it ii)] $\frac{C_{i+1}}{C_i}={\left(\frac{C}{C_i}\right)}_0$.  
\end{description} }
\end{drapeau}

 If $C$ is finite-dimensional, then by \ref{theo1}, $C$ is isomorphic to a $C_{(p_{i,j})}$ and so $C_0$ is a non-zero
 subcomodule of $C$. Moreover, if $i\geq 0$, we have
$\frac{C_{i+1}}{C_i}={(\frac{C}{C_i})}_0$, so $\frac{C_{i+1}}{C_i}$ is  non-zero   \\[1mm]
and we get in this way a flag of comodules: there is $k \in \mathbb{N}$, such that $C_0 \varsubsetneq \ldots \varsubsetneq C_k = C$. 
\newtheorem{prop1}[forme]{Proposition}
\begin{prop1}
Let $(p_{i,j})_{1\leq i \leq j \leq n}$ be a reduced family of primitive elements of type $(c_0,\ldots,c_k)$
and $(e_0,\ldots,e_n)$ the basis of $C_{(p_{i,j})}$ as decribed in  \ref{defcomo}.
Then for all  $l \in \{0 \ldots k \}$,\\ $(e_0,\ldots,e_{c_0 +\ldots +c_l-1})$ is a basis of ${(C_{(p_{i,j})})}_l$.
\end{prop1}

\noindent  {\it Proof:} as ${\cal P}=\pi_1({\cal Q}^T)$, we can write:
$$ {\cal Q}=
   \left [ \begin{tabular}{c|c|c|c}
    $Id$ & 0 & $\cdots$ & 0\\
    \hline 
    ${\cal Q}_{1,0}$ & $\ddots$ & $\ddots$ & 0\\
    \hline
    $\vdots$ & $\cdots$ & $\ddots$ &$ \ddots$\\
    \hline
    ${\cal Q}_{k,0}$ & $\cdots$ & ${\cal Q}_{k,k-1}$ & $Id$
    \end{tabular}
    \right] $$
where the diagonal blocs are in ${\cal M}_{c_0}({\cal H}_R), \ldots ,{\cal M}_{c_k}({\cal H}_R)$.
Because of coassociativity, the elements in the blocs ${\cal Q}_{i,i-1}$ are primitive, so ${\cal Q}_{i,i-1}={\cal P}_{i,i}^T$ 
and the rows of the blocs ${\cal Q}_{i,i-1}$ are linearly independent.
We easily deduce that $(e_0,\ldots,e_{c_0-1})$ is a basis of $C_0$. We conclude by induction on $n$, 
with the remark that $\frac{C}{C_0}$ is isomorphic to $C_{(p'_{i,j})}$, with:
$${\cal P}'=\left [ \begin{tabular}{c|c|c|c}
    $0$ & ${\cal P}_{2,2}$ & $\cdots$ & ${\cal P}_{2,k}$\\
    \hline
    $\vdots$ & $\ddots$ & $\ddots$ & $\vdots$\\
    \hline
    $0$ & $\cdots$ & $\ddots$ & ${\cal P}_{k,k}$\\
    \hline
    $0$ & $\cdots$ & $\cdots$ & $0$
    \end{tabular}
    \right]
    $$
 so $(p'_{i,j})$ is a reduced family of type $(c_1, \ldots, c_k)$.
\newtheorem{prop2}[forme]{Proposition}
\begin{prop2}
Let $C$ be a comodule of finite dimension with a basis $(e_0,\ldots, e_n)$ such that 
$(e_0, \ldots, e_{dim(C_i)-1})$ is a basis of $C_i$ for $0 \leq i \leq k$. Let $(p_{i,j})$ be the family of   
primitive elements built as in the proof of \ref{theo1}. Then $(p_{i,j})$ is a reduced family of type $(c_0,\ldots, c_k)$,
with  $c_0=dim(C_0)$, $c_i=dim(C_i)-dim(C_{i-1})$ for $1 \leq i \leq k$.

\end{prop2}
\noindent  {\it Proof:} as $\frac{C_i}{C_{i-1}}$ is trivial, we have:
$$ {\cal Q}=
   \left [ \begin{tabular}{c|c|c|c}
    $Id$ & 0 & $\cdots$ & 0\\
    \hline 
    ${\cal Q}_{1,0}$ & $\ddots$ & $\ddots$ & 0\\
    \hline
    $\vdots$ & $\cdots$ & $\ddots$ &$ \vdots$\\
    \hline
    ${\cal Q}_{k,0}$ & $\cdots$ & ${\cal Q}_{k,k-1}$ & $Id$
    \end{tabular}
    \right] $$
where the diagonal blocs are in ${\cal M}_{c_0}({\cal H}_R), \ldots ,{\cal M}_{c_k}({\cal H}_R)$, 
and the blocs ${\cal Q}_{i,i-1}$ are formed of primitive elements.
Suppose the rows of the bloc ${\cal Q}_{i,i-1}$ are not linearly independent.
Then we can build an element $x\in C_{i+1}-C_{i}$, with $\Delta_C(x)\equiv 1 \otimes x \left[ {\cal H}_R \otimes C_{i-1} \right]$,
hence $\overline{x}$ is a trivial element of $\frac{C}{C_{i-1}}$, which contradicts the definition of $C_i$.
We conclude using the equality ${\cal P}=\pi_1({\cal Q}^T)$.
\newtheorem{cors}[forme]{Corollary}
\begin{cors}
For any finite-dimensional comodule $C$, there exists a reduced family $(p_{i,j})$ such that $C$ is isomorphic to $C_{(p_{i,j})}$.\\ 
If $(p_{i,j})$ and $(p'_{i,j})$ are reduced families with  $C_{(p_{i,j})}$ and $C_{(p'_{i,j})}$ isomorphic, then  
$(p_{i,j})$ and $(p'_{i,j})$ have the same type.
\end{cors}
\indent In the following, we call "type of a comodule $C$" the type of any  reduced family $(p_{i,j})$ such that $C$ is isomorphic to $C_{(p_{i,j})}$.
Given $(c_0,\ldots, c_k)$, we call 
$$ G_{(c_0,\ldots, c_k)}= \left\{
\left [ \begin{tabular}{c|c|c|c}
    $g_{0,0}$ & $g_{0,1}$ & $\cdots$ & $g_{0,k}$\\ \hline
    $\vdots$ & $\ddots$ & $\ddots$ & $\vdots$\\ \hline
    0 & $\cdots$ & $\ddots$ & $g_{k-1,k}$\\ \hline
    0 & $\cdots$ & $\cdots$ & $g_{k,k}$
    \end{tabular}
    \right]
    ,g_{i,i} \in GL_{c_i}(\mathbb{Q})
    \right\}
    \subset GL_{c_0+ \ldots + c_k}(\mathbb{Q}).   
$$
$G_{(c_0,\ldots, c_k)}$ is a parabolic subgroup of $GL_{c_0+ \ldots + c_k}(\mathbb{Q})$, and it acts on the set of reduced families of type $(c_0,\ldots,c_k)$
by $g.{\cal P} = g{\cal P}g^{-1}$, where $g \in  G_{(c_0,\ldots, c_k)}$, and ${\cal P}$ is the matrix of a reduced family $(p_{i,j})$.

\newtheorem{class}[forme]{Theorem}
\begin{class}
Let $(p_{i,j})$ and $(p'_{i,j})$ be two reduced families of primitive elements of ${\cal H}_R$, 
and $(c_0,\ldots, c_k)$ be the type of $(p_{i,j})$. Then
$ C_{(p_{i,j})} \approx C_{(p'_{i,j})} $ if and only if 
$(p_{i,j}),(p'_{i,j})$ have the same type and there exists $g\in G_{(c_0,\ldots, c_k)}$, such that ${\cal P'}=g.{\cal P} $.

\end{class}

\noindent  {\it Proof:}
we put $C=C_{(p_{i,j})},C'= C_{(p'_{i,j})} $.\\
$\Leftarrow$: we have ${\cal P'}=g.{\cal P}$, so ${\cal Q}={(g^{T})}^{-1}{\cal Q'}g^T$. 
Let ${(g^{T})}^{-1}=(a_{i,j})_{0 \leq i,j \leq n}$, $g^T=(b_{i,j})_{0 \leq i,j \leq n}$ and let $(f_0, \ldots f_n)$ be the basis of $C$
defined by $f_i=\sum_j b_{i,j}e_j$.
An easy direct computation shows that $\Delta_C(f_i)=\sum_{j,k} (b_{i,j}  Q_{j,k} a_{k,l}) \otimes f_l
               = \sum_{i}  Q'_{i,l} \otimes f_l$. So $C \approx C'$.\\
\\
$\Rightarrow$: then there exists $A \in GL_{n+1}(\mathbb{Q})$, with inverse $B$ such that 
if $f_i = \sum_j b_{i,j} e_j$, then $\Delta_C(f_i)=\sum_l {\cal Q'}_{i,l} \otimes f_l$.
Then the same computation shows that $Q'_{i,l}= \sum_{j,k} b_{i,j} Q_{j,k} a_{k,l}$ or equivalently:
${\cal Q'}=A^{-1} {\cal Q} A$. Hence, ${\cal P}=A^T {\cal P'} {A^T}^{-1}$.   
As $(p'_{i,j})$ is reduced, $C_i=(f_0, \ldots, f_{c_0 + \ldots + c_i -1}) 
                                =(e_0, \ldots, e_{c_0 + \ldots + c_i -1})$
so $ A^T \in  G_{(c_0,\ldots, c_k)}$. \\

We have now entirely proved the following theorem:
\newtheorem{resume}[forme]{Theorem}
\begin{resume}
Let ${\cal P}_{(c_0 \ldots c_k)}$ be the set of the reduced families of primitive elements of ${\cal H}_R$
of type $(c_0,\ldots, c_k)$, and ${\cal O}_{(c_0,\ldots,c_k)}$ the orbit space under the action of the parabolic subgroup $G_{(c_0,\ldots, c_k)}$ of $GL_{c_0+ \ldots + c_k}(\mathbb{Q})$.
Then there is a bijection from ${\cal O}_{(c_0 \ldots c_k)}$  into the set of ${\cal H}_R$-comodules of type $(c_0, \ldots, c_k)$.
Moreover there is a bijection from the disjoint union of the ${\cal O}_{(c_0 \ldots c_k)}$'s into the set of finite-dimensional comodules. 
\end{resume}

 {\it Example:} let $C$ be a comodule of dimension 2. Then its type can be $(2)$ or $(1,1)$. We have:
$$ {\cal P}_{(2)}= \left\{ \left[ \begin{array}{cc}
                                         0&0\\ 0&0\\ \end{array} \right] \right\}, \quad
{\cal P}_{(1,1)}= \left\{ \left[ \begin{array}{cc}
                                         0& p\\ 0&0\\ \end{array} \right]/ p \neq 0 \right\}. $$
Let $\left[ \begin{array}{cc} 0& p\\ 0&0\\ \end{array} \right]$ and $\left[ \begin{array}{cc} 0& p'\\ 0&0\\ \end{array} \right] $ $\in {\cal P}_{(1,1)}$. They are in the same orbit under the action of $G_{(1,1)}$ if and only if $\exists\lambda \in \mathbb{Q}^*$, $p'=\lambda p$.
Hence, ${\cal O}_{(1,1)}$ is in bijection with the projective space associated to $Prim(\HR)$, and ${\cal O}_{(2)}$ is reduced to a single point, which corresponds to the trivial comodule of dimension 2.
\\

 We now give a caracterization of comodules of type $(n+1)$ and type $(1,\ldots,1)$.
\clearpage
\newtheorem{Types}[forme]{Proposition}
\begin{Types}
Let $C$ be a comodule of dimension $n+1$.\\[1mm]
\begin{tabular}{clcl}
1. & $C$ is of type $(n+1)$ & $\Longleftrightarrow$  & $C$ is trivial.\\[1mm]
2. & $C$ is of type $(1, \ldots,1)$ & $\Longleftrightarrow$ &
$\forall i\in \{1 \ldots n+1\}$, $C$ has a unique subcomodule\\
& & & of dimension $i$.
\end{tabular} \\
In particular, if $C$ is of type $(1, \ldots,1)$, $C$ admits a unique complete flag of subcomodules. 
\end{Types}
\noindent{\it Proof:} $1.$ is obvious.\\
2. $\Leftarrow$: let $C^{\,(i)}$ be the unique subcomodule of dimension $i+1$ of $C$. Let $x \in C_0$, $x\neq 0$. Then $(x)$ is a subcomodule of dimension $1$ of $C$, so $(x)=C^{\,(0)}$ and we get $C_0=C^{\,(0)}$.  
Suppose that $C_{i-1}=C^{\,(i-1)}$. Let $x \in C_i-C_{i-1}$, then $C_{i-1} \oplus (x)$ is a subcomodule of dimension $i+1$ of $C$, so it is equal to $C^{\,(i)}$ and we get $C_i=C^{\,(i)}$. Hence, the type of $C$ is $(1,\ldots,1)$.

 $\Rightarrow$: let $C'$ be a subcomodule of dimension $1$ of $C$. Then $C'$ is trivial, so $C' \subset C_0$. As $dim(C_0)=1$, $C'=C_0$. Suppose that $C$ has a unique subcomodule of dimension $i$. Then it is $C_{i-1}$. Let $C''$ be a subcomodule of dimension $i+1$.
It has a subcomodule of dimension $i$, so $C_{i-1} \subset C''$. Moreover, $\frac{C''}{C_{i-1}}$ is trivial, so $C''\subset C_{i}$. As they have the same dimension, $C''=C_{i}$.
\\

To conclude this section, we indicate how finite-dimensional comodules can help in renormalization.
Recall the Toy model of \cite{Connes}.
For a rooted tree $t$ with $n$ vertices, enumerated such that the root has number one, we associate the integral

$$x_t(c)= \int_0^{\infty} \frac{1}{y_1+c} \prod_{i=2}^{n} \frac{1}{y_i+y_{j(i)}} \, y_n^{-\varepsilon} d y_n \ldots y_1^{-\varepsilon} d y_1, \: \forall c>0,$$
where $j(i)$ is the number of the vertex to which the $i$-th vertex is connected via its incomming edge.
   
  Let $\{t_1, \ldots, t_m=t\}=\{R^C(t)/\mbox{$C$ cut of $t$}\}$.
  We take the comodule ${\cal C}$ with basis $(x_{t_1}, \ldots x_{t_m})$, and  structure map defined by
  $$\Delta_{\cal C}(x_{t_i}) = 1 \otimes x_{t_i} +\sum_{\mbox{admissible cuts C of $t_i$}}
  P^C(t_i) \otimes x_{R^C(t_i)}. $$
  With $[M]=x_M(0)$ for $M$ a non-empty forest, and $[1]=1$, we consider the integral:
  $$ \overline{x_t}(c)=\left( ([.] \otimes Id) \circ (S \otimes Id) \circ(\Delta_{\cal C}) \right)(x_t)$$
  Then the renormalized function is: 
  $$x_t^R(c)=
\begin{array}{c}
\textnormal{lim}\\[-2mm]
\mbox{\footnotesize{$\varepsilon \longmapsto 0$}}
\end{array}
 \: \mbox{\large{(}}\overline{x_t}(c) -[\overline{x_t}(c)]\mbox{\large{)}}.$$
We don't have anymore to worry about non commutativity within the forests. 
\\

{\it Example:}
\begin{eqnarray*}
  x_{l_1}(c)&=& \int_0^{\infty} \frac{1}{y_1+c} \: y_1^{-\varepsilon} d y_1,\\
 x_{l_2}(c)&=& \int_0^{\infty} \frac{1}{y_1+c} \:\frac{1}{y_2+y_1} \: y_2^{-\varepsilon} d y_2 \: y_1^{-\varepsilon} d y_1,\\
x_{l_3}(c)&=& \int_0^{\infty} \frac{1}{y_1+c} \:\frac{1}{y_2+y_1} \:\frac{1}{y_3+y_2} \: y_3^{-\varepsilon} d y_3 \: y_2^{-\varepsilon} d y_2 \: y_1^{-\varepsilon} d y_1.
\end{eqnarray*}

We take the comodule ${\cal C}$ with basis $(x_{l_1},x_{l_2},x_{l_3})$. We then get:
$$ \Delta_{\cal C}(x_{l_3})= 1 \otimes x_{l_3} + l_1 \otimes x_{l_2} + l_2 \otimes x_{l_1}.$$
$$\mbox{So } \overline{x_{l_3}}(c)= x_{l_3}(c) - [x_{l_1}(c)]x_{l_2}(c) - [x_{l_2}(c)]x_{l_1}(c) + [x_{l_1}(c) x_{l_1}(c)]x_{l_1}(c),$$
\\[-6mm]
$$\mbox{and } x_{l_3}^R(c) = 
\lim_{\varepsilon\rightarrow 0}
\left(
\begin{array}{ccccccc}
x_{l_3}(c)& -& [x_{l_1}(c)]x_{l_2}(c) &-& [x_{l_2}(c)]x_{l_1}(c) &+& [x_{l_1}(c) x_{l_1}(c)]x_{l_1}(c) \\[2mm]
-[x_{l_3}(c)] &+& [[x_{l_1}(c)]x_{l_2}(c)] &+& [[x_{l_2}(c)]x_{l_1}(c)] &-& [[x_{l_1}(c) x_{l_1}(c)]x_{l_1}(c)]
\end{array}
\right).$$

\section{${\cal H}_R$ as a comodule. Bigrading ${\cal H}_R$}
\esp Here, we consider the (left-)comodule $C=({\cal H}_R,\Delta)$. Of course it is not finite-dimensional, but it is the union of finite-dimensional comodules 
(for example, the comodules linearly spanned by the forests of weight less than $n,n \in \mathbb{N}$).
\newtheorem{grad}[forme]{Proposition}
\begin{grad}
$C_0=(1)$; if $i \geq 1$ then $C_i=(1) \oplus \oplus_{j=1}^{j=i} Im(F_j)$.  
\end{grad}
\noindent  {\it Proof:} $C_0$: let $x \in C, \Delta(x) = 1 \otimes x$. Then $x =(Id \otimes \varepsilon)(\Delta(x))=\varepsilon(x)1$: $x$ is cons\-tant.  
$i \geq 1$: induction on $i$. Let $x \in C_{i+1},\Delta(x)=1 \otimes x + x \otimes 1
  + \sum_j x^{(1)}_{j} \otimes x^{(2)}_j$. By hypothesis, the $x^{(2)}_j$'s are in $C_i=(1) \oplus \oplus_{j=1}^{j=i} Im(F_j)$. 
  Suppose that $x^{(2)}_1 \ldots x^{(2)}_l$ are in $Im(F_i)$, the others in $C_{i-1}$.
  By coassociativity of $\tilde{\Delta}$, $x^{(1)}_1 \ldots x^{(1)}_l$ are primitive.\\ 
  $$\hspace*{-4.2cm}\mbox{Then }\Delta \left(x - F_{i+1} \left( \sum_{j=1}^{j=l} x^{(1)}_j \otimes F_i^{-1}(x^{(2)}_j) \right) \right)  
            \equiv 1 \otimes x \left[ {\cal H}_R \otimes C_{i-1} \right],$$
        $$\mbox{so } x - F_{i+1} \left( \sum_{j=1}^{j=l} x^{(1)}_j \otimes F_i^{-1}(x^{(2)}_j) \right) \in C_i.  
        \mbox{ Hence, }C_{i+1} = C_{i} + Im(F_{i+1})\mbox{. The result is then}$$ trivial.          
\newtheorem{deco}[forme]{Proposition}
\begin{deco}
\label{dec}
$C = (1) \oplus \oplus_{j=1}^{j=\infty} Im(F_j).$
\end{deco}

\noindent  {\it Proof:} let ${\cal H}_n$ be the subspace of ${\cal H}_R $ generated by the homogeneous elements of weight $n$. 
Then $ \oplus_{i=0}^{n}{\cal H}_i$ is a subcomodule of $C$. By \ref{drap}, we have
$( \oplus_{i=0}^{n}{\cal H}_i)_k \subset C_k$. For a $k$ great enough, we have:
$ \oplus_{i=0}^{n}{\cal H}_i=( \oplus_{i=0}^{n}{\cal H}_i)_k \subset C_k$.
So as  ${\cal H}_R=\oplus_{i=0}^{\infty}{\cal H}_i$, we have the result.\\

It is now easy to see that $C_i=Ker(\tilde{\Delta}^{i})\oplus(1)$. We recognize then 
the second grading of \cite{Kreimer}, that is to say $C_i=\{x \in {\cal H}_R/deg_p(x) \leq i\}$, which defines $deg_p$.
Following \cite{Kreimer}, we put ${\cal H}_{n,k}={\cal H}_n \cap C_k$,  
  $h_{n,k}=dim({\cal H}_{n,k})$, and $r_n=dim({\cal H}_n)$. One has $h_{0,0}=1$ and $h_{n,0}=0$ if $n \neq 0$.
Note that $h_{n,1}=dim({\cal H}_n \cap Prim({\cal H}_R))$.

\newtheorem{dimension}[forme]{Proposition}
\begin{dimension}
\label{dimen}
$$ \mbox{Let }\Theta_n=\sum_{b_1 +2b_2+ \ldots +n b_n=n} (-1)^{b_1 + \ldots +b_n +1}
             \frac{(b_1 + \ldots +b_n)!}{b_1! \ldots b_n!}X_1^{b_1} \ldots X_n^{b_n} \in \mathbb{Q}[X_1 \ldots X_n ]$$
$$ \mbox{and }\varphi_{n,k}=\sum_{\begin{array}{c}
\scriptstyle{b_1 +2b_2+ \ldots +n b_n=n}\\
\scriptstyle{b_1+b_2+ \ldots +b_n=k}
\end{array} } 
\frac{k!}{b_1! \ldots b_n!}X_1^{b_1} \ldots X_n^{b_n} \in \mathbb{Q}[X_1 \ldots X_n ].$$

Then $h_{n,1}=\Theta_n(r_1,\ldots,r_n) \, \forall n \in \mathbb{N}$, and  
  $h_{n,k}=\varphi_{n,k}(h_{1,1},\ldots,h_{n,1}) \, \forall n,k \in \mathbb{N}^*$.
\end{dimension}

\noindent  {\it Proof:} 
$$ \hspace*{-1cm} \mbox{We also need }\Phi_n=\sum_{b_1 +2b_2+ \ldots +n b_n=n} 
             \frac{(b_1 + \ldots +b_n)!}{b_1! \ldots b_n!}X_1^{b_1} \ldots X_n^{b_n} \in \mathbb{Q}[X_1 \ldots X_n ].$$

As the $F_i$ are homogeneous, we have ${\cal H}_n=\oplus_{i=1}^{n} \oplus_{b_1+ \ldots +b_i=n } F_i(\otimes_{j=1}^{n} {\cal H}_{b_j,1} )$. 
As the $F_i$ are injective, we find: $r_n=\Phi_n(h_{1,1}, \ldots,h_{n,1})$. 
Let's work in the algebra of formal power series $\mathbb{Q}[[X_1, \ldots,X_n,\ldots]]$.
In this algebra,we have:
\begin{eqnarray*}
\sum_{(b_1, \ldots, b_n) \neq ( 0, \ldots,0)}\frac{(b_1 + \ldots +b_n)!}{b_1! \ldots b_n!}X_1^{b_1} \ldots X_n^{b_n}
&=&\sum_{k \neq 0} \left(\sum_{b_1 + \ldots + kb_k = k}\frac{(b_1 + \ldots +b_k)!}{b_1! \ldots b_k!}X_1^{b_1} \ldots X_k^{b_k}\right)\\
& =&\sum_{k \neq 0} \Phi_k(X_1,\ldots,X_k)\\
&=&\sum_{l \neq 0} \left(\sum_{b_1 + \ldots + b_k = l}\frac{l!}{b_1! \ldots b_k!}X_1^{b_1} \ldots X_k^{b_k} \right)\\
&=&\sum_{l \neq 0} \left( \sum_{i \neq 0} X_i \right)^l
= \frac{ \sum_{i \neq 0} X_i }{1-\sum_{i \neq 0} X_i}.
\end{eqnarray*}
We then get:
$$ \sum_{k \neq 0} \Phi_k(-\Phi_1, \ldots, -\Phi_k)=
     \frac{- \sum_{i \neq 0} \Phi_i}{1 + \sum_{i \neq 0} \Phi_i}
     = \frac{- \frac{ \sum X_i }{1-\sum X_i}}{1 + \frac{ \sum X_i }{1-\sum X_i}}    
     =-\sum_{i \neq 0} X_i.$$

 Hence, by putting $X_i$ in weight $i$ and by comparing the homogeneous parts of each member, 
we find $\Phi_k(-\Phi_1, \ldots, -\Phi_k)=-X_k$, or equivalently $\Theta_k(\Phi_1, \ldots, \Phi_k)=X_k$.
So $\Theta_k\left( \Phi_1(h_{1,1}), \ldots, \Phi_k(h_{1,1}, \ldots,h_{k,1}) \right)=
       \Theta_k(r_1, \ldots,r_k)=h_{k,1}$.\\
If $ k >1 $, then ${\cal H}_{n,k}= \oplus_{c_1+\ldots +c_k=n } 
       F_k({\cal H}_{c_1,1}, \ldots, {\cal H}_{c_k,1})$.
As $F_k$ is injective, we find the announced result.\\

 We denote $H(X,Y)=\sum_{n,k} h_{n,k}X^nY^k$, $H_j(x)=\sum_n h_{n,j}X^n$, $R(X)=\sum_n r_nX^n$. 
\\
The second formula of \ref{dimen} implies that $H_j(X)=H_1(X)^j, \forall j \in \mathbb{N}$.
The first formula implies that $1-H_1(X)=\frac{1}{R(X)}$.
We have then $H(X,Y)=\sum_{j=0}^{\infty} H_j(X)Y^j= \sum_{j=0}^{\infty} 
\left[ H_1(X)Y \right]^j=\frac{1}{1-H_1(X)Y}=\frac{R(X)}{Y+(1-Y)R(X)}$, which is a reformulation of the main
theorem of \cite{Kreimer} (with a small difference because of the different definitions of $R(X)$). We give the first values of $r_n$ and $h_{n,1}$ in the appendix (see also \cite{Sloane}).

\section{The Lie algebra ${\cal L}^1$}

\newtheorem{Ulibre}[forme]{Proposition}
\begin{Ulibre}
\label{Ulib}
\begin{enumerate}
\item $\Ul$ is a free algebra;
\item $ \forall l_1,l_2\in\Ul$, $weight(l_1l_2)=weight(l_1) +weight(l_2)$.
\end{enumerate}
\end{Ulibre}
{\it proof:}
let $(p_i)_{i \geq 1}$ be a basis of $Prim(\HR)$ such that the $p_i$'s are homogeneous for the weight.
By proposition \ref{dec} and lemma \ref{Fis}, $(p_{i_1} \top \ldots \top p_{i_k})_{k \geq 0, i_1, \ldots , i_k \geq 1}$ is a basis of $\HR$.
We define $f_{j_1, \ldots ,j_l} \in \HR^*$ by :
$$ f_{j_1, \ldots ,j_l}(p_{i_1} \top \ldots \top p_{i_k})=\left\{
\begin{tabular}{cl}
1& if  $(j_1, \ldots, j_l)= (i_1,\ldots,i_k) $ \\
0& if $ (j_1, \ldots, j_l)\neq (i_1,\ldots,i_k). $
\end{tabular}
\right.$$
As the $(p_{i_1} \top \ldots \top p_{i_k})$'s are homogeneous for the weight, $(f_{j_1, \ldots ,j_l})_ {k \geq 0, i_1, \ldots , i_k \geq 1}$ is a basis of $\HR^{*g}$.
\begin{eqnarray*}
 (f_{j_1,\ldots,j_l}f_{j'_1,\ldots,j'_n}, p_{i_1}\top \ldots \top p_{i_k})&=&(f_{j_1,\ldots,j_l} \otimes f_{j'_1,\ldots,j'_n}, \Delta(p_{i_1}\top \ldots \top p_{i_k})) \\
&=&(f_{j_1,\ldots,j_l} \otimes f_{j'_1,\ldots,j'_n}, \sum_{s=0}^{k}p_{i_1}\top \ldots \top p_{i_s}\otimes p_{i_{s+1}}\top \ldots \top p_{i_k})\\
&=& \left\{ \begin{array}{cl}
            1& \mbox{ if }(j_1,\ldots, j_l,j'_1,\ldots,j'_n)=(i_1,\ldots,i_k),\\
            0& \mbox{ if }(j_1,\ldots, j_l,j'_1,\ldots,j'_n)\neq (i_1,\ldots,i_k).\\
      \end{array}
\right.
\end{eqnarray*}
So $f_{j_1,\ldots,j_l}f_{j'_1,\ldots,j'_n}=f_{j_1,\ldots,j_l,j'_1,\ldots,j'_n}$, hence $\HR^{*g}$ and the free algebra generated by the $f_i$'s, $i\geq 1$, are isomorphic algebras.  
Moreover, the $f_i$'s are homogeneous elements of $\HR^{*g}$, so we have $weight(ff')=weight(f)+weight(f')$ $ \forall f,f' \in \HR^{*g}$. As $\HR^{*g}$ and $\Ul$ are isomorphic graded algebras, the proposition is proved.\\

Let ${\cal A}$ be the augmentation ideal of $\Ul$, that is to say ${\cal A}=ker(\varepsilon)$.

\newtheorem{orthogonal}[forme]{Lemma}
\begin{orthogonal}
\label{ortho} in the duality between $\Ul$ and $\HR$, the orthogonal of ${\cal A}^2 \oplus (1)$ is $Prim(\HR)$.
\end{orthogonal}
{\it Proof:}\\
$({\cal A}^2\oplus(1))^{\perp} \subseteq Prim(\HR)$: let $x \in  ({\cal A}^2\oplus(1))^{\perp}$, and $l_1,l_2\in \Ul$. One has to show that $(l_1 \otimes l_2,\Delta(x))=(l_1 \otimes l_2,x\otimes 1 + 1\otimes x)$, that is to say
$(l_1l_2,x)=\varepsilon(l_1)(l_2,x) + (l_1,x)\varepsilon (l_2)$. As $\Ul=(1) \oplus ker(\varepsilon)$, one has four cases to considerate:
\begin{enumerate}
\item $\varepsilon(l_1)=\varepsilon(l_2)=0$: then $l_1l_2 \in {\cal A}^2$, so $(l_1l_2,x)=0=\varepsilon(l_1)(l_2,x) + (l_1,x)\varepsilon (l_2)$;
\item $\varepsilon(l_1)=0$ and $l_2=1$: obvious;
\item $l_1=1$ and $\varepsilon(l_2)=0$: obvious;
\item $l_1=l_2=1:$ one has to show that $(1,x)=2(1,x)$; as $(1,x)=0$ it is true.
\end{enumerate}
$Prim(\HR)\subseteq ({\cal A}^2\oplus(1))^{\perp} $: equivalently we show ${\cal A}^2\oplus(1) \subseteq (Prim(\HR))^{\perp}$.
Let $p \in Prim(\HR)$, then $(1,p)=\varepsilon(p)=0$, so $1 \in (Prim(\HR))^{\perp}$. Let $l \in {\cal A}^2$. One can suppose that $l=l_1l_2,$ $\varepsilon(l_1)=\varepsilon(l_2)=0$.
Let $p \in Prim(\HR)$. 
$$ (l_1l_2,p)=(l_1 \otimes l_2, p\otimes 1 + 1 \otimes p)= \varepsilon(l_1)(l_2,p)+(l_1,p)\varepsilon(l_2)=0.$$         

We denote by $\Ul_n$ the space of the homogeneous elements of $\Ul$ of weight $n$. We have $dim(\Ul_n)=dim({\cal H}_n)=r_n$. Moreover, ${\cal A}=\oplus_{n\geq 1} \Ul_n$.
We denote ${\cal A}^2_n={\cal A}^2 \cap \Ul_n$. We have ${\cal A}^2 = \oplus_{n\geq 1} {\cal A}^2_n$. Now, observe that ${\cal L}^1 + {\cal A}^2=\Ul$
(it is obvious when one takes a Poincar\'e-Birkhoff-Witt basis of $\Ul$).
So for every $n \geq 1$, we can choose a subspace $G_n$ of ${\cal L}^1$, such that $\Ul_n=G_n \oplus {\cal A}^2_n$.
By lemma \ref{ortho}, $dim(G_n)=dim(Prim(\HR)\cap {\cal H}_n)=h_{n,1}$. We denote $G=\oplus_{n\geq 1}G_n$.

\newtheorem{genL1}[forme]{Lemma}
\begin{genL1}
$G$ generates the algebra $\Ul$.
\end{genL1}
{\it Proof:}
we denote by $<G>$ the subalgebra of $\Ul$ generated by $G$. Let $l \in \Ul$, homogeneous of weight $n$; we proceed by induction on $n$.
If $n=0$, then $l$ is constant: it is then obvious. Suppose that every element of weight less than $n$ is in $<G>$.
As $\Ul_n=G_n \oplus {\cal A}^2_n$, one can suppose that $l=l_1l_2$, with $l_1, l_2\in {\cal A}$. By lemma  \ref{Ulib},
$weight(l)=weight(l_1)+weight(l_2)$, so $weight(l_1)<n$, and $weight(l_2)<n$. Then $l_1, l_2 \in <G>$, and $l \in <G>$.\\

We denote by ${\cal F}(G)$ the free associative algebra generated by the space $G$. The gradation of $G$ induces a gradation of the algebra ${\cal F}(G)$.
By the last lemma, we have a surjective algebra morphism:
$$\Upsilon:
\left\{  \begin{array}{ccc}
          {\cal F}(G)& \longmapsto & \Ul \\
           g \in G & \longmapsto & g.\\
\end{array} \right.
$$
Moreover, $\Upsilon$ is homogeneous of degree 0. We now calculate the dimension $f_n$ of the homogeneous part of weight $n$ of ${\cal F}(G):$

$$ \hspace{-2.7cm} f_n= \sum_{\begin{array}{c}
\scriptstyle{               a_1+\ldots+a_k=n}\\
\scriptstyle{             a_i \geq 1 \: \forall i} \end{array}  }
    h_{a_1,1} \ldots h_{a_k,1}$$
$$= \sum_{b_1+2b_2+\ldots+nb_n=n} \frac{b_1+ \ldots +b_n}{b_1! \ldots b_n!} h_{1,1}^{b_1}\ldots h_{n,1}^{b_n}= r_n.
$$
(For the second equality, $b_i$ is the number of the $a_j$'s equal to $i$; the third equality was shown in the proof of proposition \ref{dimen}).\\
As the homogeneous parts of $\Ul$ and ${\cal F}(G)$ have the same finite dimension, and as $\Upsilon$ is surjective and homogeneous of degree $0$, it is in fact an isomorphism.\\

We now put a Hopf algebra structure on ${\cal F}(G)$ by putting $\Delta(g)=g \otimes 1 + 1 \otimes g$ $\forall g \in G$. As $G \subset {\cal L}^1$, the elements of $G$ are primitive in both $\Ul$ and ${\cal F}(G)$, so $\Upsilon$ is a Hopf algebra isomorphism.
Hence, it induces a Lie isomorphism between $Prim({\cal F}(G))$ and $Prim(\Ul)={\cal L}^1$. But  $Prim({\cal F}(G))$ is isomorphic to the free Lie algebra generated by $G$ (see for example \cite{Bourbaki}), so we have the following result:

\newtheorem{Llibre}[forme]{Theorem}
\begin{Llibre}
${\cal L}^1$ is a free Lie algebra.
\end{Llibre}

\section{Primitive elements}
\subsection{Primitive elements of the Hopf Algebra ${\cal H}_{ladder}$}

\esp First, we construct a family of primitive elements of ${\cal H}_{ladder}$. 
For that, we introduce the Hopf algebra $\mathbb{Q}[X_1,\ldots,X_n, \ldots]$ with coproduct 
defined by $ \Delta(X_i)=X_i \otimes 1 + 1 \otimes X_i$. In this algebra let 
$$\Psi_n=\sum_{a_1+2 a_2\ldots +na_n=n} \frac{X_1^{a_1} \ldots X_n^{a_n}}{a_1! \ldots a_n!}.$$
\newtheorem{lemme3}[forme]{Lemma}

\begin{lemme3}
\label{lem3}
$ \Delta(\Psi_n)=\sum_{j=0}^{j=n} \Psi_j \otimes \Psi_{n-j}$.       
\end{lemme3}
\noindent  {\it Proof:}
one easily shows that: $$\Delta(X_1^{a_1} \ldots X_n^{a_n})=\sum_{i=0}^{n} \sum_{k_i=0}^{k_i=a_i} \binom{a_1}{k_1} \ldots \binom{a_n}{k_n} X_1^{k_1} \ldots X_n^{k_n} \otimes X_1^{a_1-k_1} \ldots X_n^{a_n-k_n}. $$
So 
\begin{eqnarray*}
 \Delta(\Psi_n)&=& \sum_{a_1+2 a_2\ldots +na_n=n} \: \sum_{i=0}^{n}\: \sum_{k_i=0}^{k_i=a_i} \frac{1}{a_1! \ldots a_n!} \binom{a_1}{k_1} \ldots \binom{a_n}{k_n} X_1^{k_1} \ldots X_n^{k_n} \otimes X_1^{a_1-k_1} \ldots X_n^{a_n-k_n} \\
&=& \sum_{b_1 + \ldots +nb_n + c_1+ \ldots +nc_n=n} \frac{\binom{b_1+c_1}{b_1} \ldots \binom{b_n+c_n}{b_n}}{(b_1+c_1)! \ldots (b_n+c_n)!} X_1^{b_1} \ldots X_n^{b_n} \otimes X_1^{c_1} \ldots X_n^{c_n}\\
&=& \sum_{j=0}^{n} \esp  \sum_{b_1 + \ldots +jb_j=j} \esp \sum_{c_1+ \ldots +(n-j)c_{n-j}=n-j} \frac{1}{b_1!c_1! \ldots b_n!c_n!} X_1^{b_1} \ldots X_j^{b_j} \otimes X_1^{c_1} \ldots X_{n-j}^{c_{n-j}}\\
&=& \sum_{j=0}^n \Psi_j \otimes \Psi_{n-j}.
\end{eqnarray*}

We define a sequence $(P_i)_{i \geq 1}$ of elements in  ${\cal H}_{ladder}$ by:
$$ P_1=l_1,\: P_n=l_n-\Psi_n(P_1,\ldots,P_{n-1},0) \mbox{ } \forall n \geq 2.$$

 As $\Psi_n=X_n + \Psi_n(X_1,\ldots,X_{n-1},0)$, we have $l_n= \Psi_n(P_1,\ldots,P_{n-1},P_n)$.
\newtheorem{suitepi}[forme]{Proposition:}
\begin{suitepi}
$P_i$ is primitive for all $i \geq 1$.
\end{suitepi}

\noindent  {\it Proof:} induction on $i$.
It is trivial for $i=1$. Suppose it is true for each $j \leq i-1$.
Then 
\begin{eqnarray*}
 \tilde{\Delta}(l_i)&=&\sum_{j=1}^{j=i-1} l_j \otimes l_{i-j}\\
& =&\sum_{j=1}^{j=i-1} \Psi_j(P_1,\ldots,P_j) \otimes \Psi_{i-j}(P_1,\ldots,P_{i-j}) \\
 &=&\tilde{\Delta}\left(\Psi_i (P_1, \ldots,P_{i-1},0)\right) 
\end{eqnarray*}
    by \ref{lem3}, and the fact that $P_1, \ldots, P_{i-1}$ are primitive.
So $\tilde{\Delta}\left(l_i -\Psi_i(P_1, \ldots,P_{i-1},0)\right)=\tilde{\Delta}(P_i)=0$, hence 
$P_i$ is primitive.\\

 We work again in $\mathbb{Q}[[X_1, \ldots,X_n,\ldots]]$. In this algebra, we have:
\begin{eqnarray*}
\sum_{(b_1, \ldots, b_n) \neq ( 0, \ldots,0)}\frac{X_1^{b_1} \ldots X_n^{b_n}}{b_1! \ldots b_n!}
&=&\sum_{k \neq 0} \left(\sum_{b_1 + 2b_2+\ldots + kb_k = k}\frac{X_1^{b_1} \ldots X_k^{b_k}}{b_1! \ldots b_k!}\right)\\
&=&\sum_{k \neq 0} \Psi_k(X_1,\ldots, X_k)\\
&=&\sum_{l \neq 0} \frac{1}{l!} \left(\sum_{b_1+b_2+\ldots +b_k=l} \frac{l!}{b_1! \ldots b_k!} X_1^{b_1} \ldots X_k^{b_k}\right)\\
&=&\sum_{l \neq 0} \frac{1}{l!}\left( \sum_{i \neq 0} X_i \right)^l \\
&=& (\exp-1)\left( \sum_{i \neq 0} X_i \right).
\end{eqnarray*}

So $\ln\left(1+\sum_{k \neq 0} \Psi_k(X_1,\ldots, X_k) \right)
=\ln \left(1+(\exp-1)\left( \sum_{i \neq 0} X_i \right) \right)= \left( \sum_{i \neq 0} X_i \right).$
By putting $X_i$ in weight $i$, and comparing the homogeneous parts, we find:
$$\sum_{a_1 + \ldots +ia_i=i} (-1)^{a_1+\ldots+a_i+1}
                   \frac{(a_1+\ldots +a_i-1)!}{a_1! \ldots a_i!}
                   \Psi_1^{a_1} \ldots \Psi_i^{a_i} = X_i.$$

 As $\Psi_i(P_1,\ldots,P_i)=l_i$, we deduce:

\newtheorem{dirpi}[forme]{Proposition}
\begin{dirpi}
$$P_i=\sum_{a_1 + \ldots +ia_i=i} (-1)^{a_1+\ldots+a_i+1}
                   \frac{(a_1+\ldots +a_i-1)!}{a_1! \ldots a_i!}
                  l_1^{a_1} \ldots l_i^{a_i}. $$
\end{dirpi}

In ${\cal H}_R$, consider the projection $\pi_c$ on the space spanned by rooted trees, which vanishes on  
the space spanned by non connected forests. We have:
\newtheorem{lemme4}[forme]{Lemma}
\begin{lemme4}
\label{lem4}
let $p \in {\cal H}_R$ be a primitive element such that $\pi_c(p)=0$. Then $p=0$.
\end{lemme4}
\noindent  {\it Proof:}
suppose $p \neq 0$, and let write $p=\sum_{\alpha=(\alpha_1,\ldots,\alpha_k)} 
          a_{\alpha} t_1^{\alpha_1} \ldots  t_k^{\alpha_k}$, where the $t_i$'s are rooted trees, with $\frac{\partial p}{\partial t_i} \neq 0$. One can suppose that $weight(t_k) \geq weight(t_i)$ $\forall i$.
Let $t_1^{\alpha_1}\ldots t_k^{\alpha_k}$ such that $\alpha_k \neq 0$ and $a_{\alpha}\neq 0$.

 Let $F$ a forest such that in the basis $(F_1 \otimes F_2)_{F_i \: forest}$ of ${\cal H}_R \otimes {\cal H}_R$, 
the coefficient of $t_1^{\alpha_1}\ldots t_{k_1}^{\alpha_{k-1}}\otimes t_k$ in $\Delta(F)$
is $\neq 0$. Then $F=t_1^{\alpha_1} \ldots  t_k^{\alpha_k}$, and then this coefficient is $\alpha_k$, or there exists  $t'$ a rooted tree with $weight(t')>weight(t)$, such that $\frac{\partial F}{\partial t'} \neq 0$.
So the coefficient of $t_1^{\alpha_1}\ldots t_{k_1}^{\alpha_{k-1}}\otimes t_k$ in $\Delta(p)$ is $\alpha_k a_{\alpha} \neq 0$.
As $p$ is primitive, $t_k=1$ or $t_1^{\alpha_1}\ldots t_{k_1}^{\alpha_{k-1}}=1$. If $t_k=1$, then $p$ is constant: this is a contradiction, because then $p$ would not be primitive.
So $t_1^{\alpha_1}\ldots t_{k_1}^{\alpha_{k-1}}=1$, and then $\pi_c(p) \neq 0$.\\

\newtheorem{basepl}[forme]{Theorem}
\begin{basepl}
$(P_i)_{i \in \mathbb{N}^*}$ is a basis of the space of primitive elements in ${\cal H}_{ladder}$.
\end{basepl}

\noindent  {\it Proof:} let p be a primitive element in ${\cal H}_{ladder}$. Then $\pi_c(p)$ is a linear combination 
of ladders, so there is a linear combination $p'$ of $P_i$ such that $\pi_c(p)=\pi_c(p')$.
By the lemma, $p=p'$.

\subsection{The operator $\pi_1$}
\esp Recall that $\pi_1$ is the projection on $Im(F_1)=Prim(\HR)$ which vanishes on $(1) \oplus \oplus_{j \geq 2} Im(F_j)$.

\newtheorem{project}[forme]{Theorem}
\begin{project}
\label{proj}
Let $F$ be a non-empty forest. $$ \mbox{We put }\tilde{\Delta}(F)=\sum_{(F)}F^{(1)}\otimes F^{(2)}\mbox{; then:}$$
$$ \pi_1(F)=F-\sum_{(F)} F^{(1)} \top \pi_1(F^{(2)}).$$
\end{project}

\noindent {\it Proof:} induction on $weight(F)$. If $weight(F)=1$, it is obvious. Suppose the formula is true for every forests of weight less than or equal to $n-1$.
Let $F$ be a forest of weight $n$. Then $weight(F^{(2)})<weight(F)$, so:
\begin{eqnarray*}
 \tilde{\Delta}(F)&=& \sum_{(F)} F^{(1)} \otimes F^{(2)}\\
&=& \sum_{(F)} F^{(1)} \otimes \left( \pi_1(F^{(2)})+\sum_{(F^{(2)})} {(F^{(2)})}^{(1)} \top \pi_1\left(({F^{(2)})}^{(2)}\right) \right)\\
& =& \sum_{(F)} \left( F^{(1)} \otimes \pi_1(F^{(2)})+ \sum_{(F^{(1)})} {(F^{(1)})}^{(1)} \otimes \left[{(F^{(1)})}^{(2)} \top \pi_1(F^{(2)})\right] \right) \mbox{ (by coassociativity)}\\
&=& \sum_{(F)} \tilde{\Delta}\left(F^{(1)} \top \pi_1(F^{(2)})\right) \mbox{ (by \ref{greffe})}.
\end{eqnarray*}
$$\mbox{So } F- \sum_{(F)} F^{(1)} \top \pi_1(F^{(2)}) \in Im(F_1) \mbox{; as }\sum_{(F)} F^{(1)} \top \pi_1(F^{(2)}) \in \oplus_{j \geq 2} Im(F_j)\mbox{, we have}$$ the result for $F$.
\\

So we have an easy way to find a family who generates the space of primitive elements of weight $n$, by induction on $n$.
Moreover, we have relations between the $\pi_1(F)$, which are given by $\pi_1(F' \top p)=0 $ for any non-empty forest $F'$  and for any  primitive element $p$ we have ever found.
So we easily have a basis  of the space of homogeneous primitive elements of weight $n$.\\

 For example, for $n=1$, we have $\pi_1(l_1)=l_1$; the basis is $(l_1)$; we have the relation $\pi_1(F'\top l_1)=0 \: \forall F'$ non-empty forest; so ${\cal R}_1:$ $\pi_1(T)=0 \: \forall T$ rooted tree of weight greater than or equal to 2.

 Hence, for $n=2$, we only have to compute $\pi_1(l_1^2)=l_1^2-2 l_1 \top \pi_1(l_1)= l_1^2-2l_2$. The basis is $(l_1^2-2l_2)$, and we have: $\pi_1(F' \top (l_1^2-l_2))=0$, which gives: ${\cal R}_2$: $\pi_1(l_1T)=0 \: \forall T$ rooted tree of weight greater than or equal to 2.

 For $n=3$, we have to compute $\pi_1(l_1^3)$; the others are zero by ${\cal R}_1$ and ${\cal R}_2$. One finds the basis $(l_1^3-3l_1l_2 +3l_3)$ and the relation ${\cal R}_3$: $\pi_1(l_1^2T)=\pi_1(l_2T) \: \forall T$ rooted tree of weight greater than or equal to 2. 

 For $n=4$, one would have to compute $\pi_1(l_1^4)$ and $\pi_1(l_2^2)$, and so on.
\\

\indent {\it Remark:} by linearity, the formula of \ref{proj} is true for any $x \in \HR$. For example, for $x=p_1p_2$, with $p_1,p_2$ primitive elements of $\HR$, one finds:
$ \pi_1(x) = p_1p_2-p_1 \top p_2 -p_2\top p_1$; hence, $ S_1(p_1)= \pi_1(-Y(p_1)l_1)$ 
         with $Y(F) = weight(F) \: F$ for all  forest $F$,
       and $S_1$ defined in \cite{Kreimer}.

\section{Classification of the Hopf algebra endomorphisms of $\HR$}
In the sequel, we will denote by $CT$ the set of (connected) rooted trees. 
\newtheorem{bla}[forme]{Definition}
\begin{bla}
\label{blabla}
Let $(P_t)_{t \in CT}$ be a family of primitive elements of $\HR$ indexed by $CT$.
Let $\Phi_{(P_t)}$ be the algebra endomorphism of $\HR$ defined by induction on $weight(T)$ by:
$$ \Phi_{(P_t)}(l_1)=P_{l_1}; $$
$$ \forall T \in CT \mbox{, with } \tilde{\Delta}(T)=\sum_{(T)} T^{\,(1)} \otimes T^{\,(2)},$$
$$\Phi_{(P_t)}(T)=\left(\sum_{(T)} \Phi_{(P_t)}(T^{\,(1)}) \top P_{T^{\,(2)}}\right)+P_T. $$
Then $\Phi_{(P_t)}$ is a bialgebra endomorphism of $\HR$.
\end{bla}
\noindent  {\it Proof:} one has to show $(\Phi_{(P_t)} \otimes \Phi_{(P_t)}) \circ \tilde{\Delta}(T)
=\tilde{\Delta}(\Phi_{(P_t)}(T)) \: \forall T \in CT$.
We proceed by induction on $n=weight(T)$. It is obvious for $n=1$, since then $T=l_1$ is primitive.
Suppose it is true for all rooted trees of weight $< n$. 
Then as $\Phi_{(P_t)}$ is an algebra endomorphism, it is true for all non connected forests of weight $\leq n$.
Let $T$ be a rooted tree of weight $n$. Then: 
\begin{eqnarray*}
\tilde{\Delta}(\Phi_{(P_t)}(T))&=&\sum_{(T)} \tilde{\Delta}\left(\Phi_{(P_t)}(T^{\, (1)}) \top P_{T^{ \, (2)}}\right)\\
& =&\left(\sum_{(T)} \Phi_{(P_t)}(T^{\, (1)}) \otimes P_{T^{ \, (2)}}\right)+ \sum_{(T)} \Phi_{(P_t)}(T^{\, (1)}) \otimes (\Phi_{(P_t)}(T^{\, (2)}) \top P_{T^{\, (3)}})\\
&=&\sum_{(T)} \Phi_{(P_t)}(T^{\, (1)}) \otimes \left[ \sum_{(T^{(2)})}\left(\Phi_{(P_t)}{((T^{\, (2)})}^{(1)}) \top P_{{(T^{\, (2)})}^{(2)}}\right) + P_T^{\, (2)} \right]\\
& =&\sum_{(T)} \Phi_{(P_t)}(T^{\, (1)}) \otimes \Phi_{(P_t)}(T^{\, (2)}).
\end{eqnarray*}

We used the induction hypothesis and \ref{greffe} for the second equality, and coassociativity of $\tilde{\Delta}$ for the third.

\newtheorem{classend}[forme]{Theorem}
\begin{classend}
Let $\Psi $ be an endomorphism of the bialgebra $\HR$. Then there exists a unique family $(P_t)$ of primitive elements, 
such that $\Psi=\Phi_{(P_t)}$.
\end{classend}

\noindent  {\it Proof:} one remarks that if $(P_t)$ and $(Q_t)$ are two families of primitive elements, such that $P_t=Q_t$ if $weight(t) \leq n$, 
then $\Phi_{(P_t)}(x) = \Phi_{(Q_t)}(x) \, \mbox{ for all  }x \mbox{ of weight } \leq n$. 
So we only have to show that there exists a family $(P_t)$ such that if we denote: 
$$ P_t^{(n)} = \left\{
 \begin{array}{ccc}
  P_t & \mbox{ if } & weight(T) \leq n \\ 
  0 & \mbox{ if } & weight(T) > n,
  \end{array} 
   \right.$$
  then $\Psi(x) = \Phi_{(P_t^{(n)})}(x)  \mbox{ for all }x \mbox{ of weight } \leq n$.
We take $P_{l_1}= \Psi(l_1)$, and then it is true for $n=1$.
Suppose we have $P_t$ for all $t$ of weight $< n$.
We put $\Phi_{(P_t^{(n-1)})}=\Phi_{n-1}$.
Let $T$ be a rooted tree of weight $n$.
\begin{eqnarray*}
 \tilde{\Delta}(\Psi(T))&=&\sum_{(T)} \Psi(T^{\,(1)}) \otimes \Psi(T^{ \, (2)})\\
& =& \sum_{(T)} \Phi_{n-1}(T^{\,(1)}) \otimes \Phi_{n-1}(T^{ \, (2)})=\tilde{\Delta}(\Phi_{n-1}(T)).
\end{eqnarray*}
We take $P_T=\Psi(T)-\Phi_{n-1}(T)$; then $\Psi(T) = \Phi_{(P_t^{(n)})}(T)$.  
\\
For the uniqueness of the family $(P_t)$, we have $\pi_1( \Psi(T))=P_T, \quad \forall T \mbox{ rooted tree}$. 
\newtheorem{bihopf}[forme]{Proposition}
\begin{bihopf}
\label{biahopf}
Let $\Psi$ be an endomorphism of the bialgebra $\HR$; then $\Psi$ is an endomorphism of the Hopf algebra $\HR$,
that is to say $\Psi \circ S =S \circ \Psi$.
\end{bihopf}
\newtheorem{Sgreffe}[forme]{Lemma}
\begin{Sgreffe}
Let $p$ be a primitive element of $\HR$ and let $x \in \HR$, with $\varepsilon(x)=0$. 
Then $$S(x \top p)= - x \top p - S(x)p - \sum_{(x)} S(x^{(1)}) (x^{(2)} \top p)  
\mbox{ where }\tilde{\Delta}(x)=\sum_{(x)} x^{(1)} \otimes x^{(2)}.$$

In particular, for $p=l_1$, $$S(B^+(x))=-B^+(x) -S(x)l_1 - \sum_{(x)} S(x^{(1)}) B^+(x^{(2)}).$$
\end{Sgreffe}
\noindent  {\it Proof:} we have $ (S \otimes Id) \circ \Delta (x) = 0$. Then we use \ref{greffe} to conclude.\\
\\
{\it Proof of the proposition:} let $F$ be a forest in $\HR$.
\begin{eqnarray*}
  \Psi \circ S(B^+(F))&=&-\Psi(F) \top P_{l_1} - \sum_{(F)} \Psi(F^{(1)}) \top P_{B^+(F^{(2)})}-P_{B^+(F)} \\
&&- \sum_{(F)} \Psi \circ S (F^{(1)}) (\Psi(F^{(2)}) \top P_{l_1})
- \sum_{(F)} \Psi \circ S(F^{(1)}) [\Psi(F^{(2)}) \top P_{B^+(F^{(3)})}]\\
&& -\Psi \circ S(F)P_{l_1}-\sum_{(F)} \Psi \circ S(F^{(1)}) P_{B^+(F^{(2)})};\\
 S \circ \Psi(B^+(F))&=&S(\Psi(F) \top P_{l_1} ) + \sum_{(F)} S( \Psi(F^{(1)}) \top P_{B^+(F^{(2)})}) + S(P_{B^+(F)})\\
&=& -\Psi(F) \top P_{l_1} -\sum_{(F)} \Psi(F^{(1)}) \top P_{B^+(F^{(2)})}- P_{B^+(F)}\\
&&- \sum_{(F)} S \circ \Psi(F^{(1)}) ( \psi(F^{(2)}) \top P_{l_1}) - \sum_{(F)} S \circ \Psi(F^{(1)}) [ \Psi(F^{(2)}) \top P_{B^+(F^{(3)})}]\\
&& - S \circ \Psi(F) P_{l_1} - \sum_{(F)} S \circ \Psi(F^{(1)}) P_{B^+(F^{(2)})}. 
\end{eqnarray*}
We conclude by an induction on the weight.

\section{Associated graded algebra of $\HR$ and coalgebra endomorphisms}

\esp As it is shown in \cite{Kreimer}, $\HR$ is filtered as Hopf algebra by $deg_p$. What is the  associated graded algebra ?

 The filtration is given by $(\HR)_n^{(P)} = \{ x \in \HR, deg_p x \leq n \}= (1) \oplus \oplus_{1}^{n} Im(F_j) =C_n =  Ker(\tilde{\Delta}^{(n)})\oplus(1)$. 
We put $\pi_i$ the projection on $Im(F_i)$ which vanishes on $(1) \oplus \oplus_{j \neq i} Im(F_j)$.

\newtheorem{termedom}[forme]{Lemma}
\begin{termedom}
\label{termedo}
Let $p_1, \ldots, p_j,p_{j+1}, \ldots, p_{j+l}$ be primitive elements of $\HR$.
Then
$$ \pi_{j+l}(p_{j+l} \top \ldots \top p_{j+1} . p_{j} \top \ldots \top p_{1})= \sum_{ \sigma \mbox{ (j,l)-shuffle}} p_{\sigma(j+l)} \top \ldots \top p_{\sigma(1)},$$
where a (j,l)-shuffle is a permutation $\sigma$ of $\{1,\ldots, j+l \: \}$ such that $\sigma(1)<\sigma(2)<\ldots<\sigma(j)$ and $\sigma(j+1)<\sigma(j+2)<\ldots<\sigma(j+l)$.

\end{termedom}

\noindent  {\it Proof:} by induction we prove:
\begin{eqnarray*}
\tilde{\Delta}^{j-l-1}(p_{j+l} \top \ldots \top p_{j+1} . p_{j} \top \ldots \top p_{1})&=& \sum_{ \sigma \mbox{ (j,l)-shuffle}} p_{\sigma(j+l)} \otimes \ldots \otimes p_{\sigma(1)}\\
& =&\tilde{\Delta}^{j-l-1}\left(\sum_{ \sigma \mbox{ (j,l)-shuffle }} p_{\sigma(j+l)} \top \ldots \top p_{\sigma(1)}\right).
\end{eqnarray*}
$$ \mbox{So } p_{j+l} \top \ldots \top p_{j+1} . p_{j} \top \ldots \top p_{1} -\sum_{ \sigma \mbox{ (j,l)-shuffle}} p_{\sigma(j+l)} \top \ldots \top p_{\sigma(1)} \mbox{ is in $(\HR)^{(P)}_{j+l-1}$,}$$ which proves the lemma.
\\

We naturally identify $(\HR)_n^{(P)} / (\HR)^{(P)}_{n-1} $ with $Im(F_n)$. We can now describe $gr(\HR)$, the associated graded Hopf algebra:\\

\begin{description}
\item[\textnormal{i)}] as vector space, $gr(\HR) = (1) \oplus \oplus_{1}^{\infty} Im(F_i);$
\item[\textnormal{ii)}] $\forall p_j \top \ldots \top p_1 \in Im(F_j),\: p_{j+l} \top \ldots \top p_{j+1} \in Im(F_l)$,\\
$$(p_{j+l} \top \ldots \top p_{j+1}) * (p_{j} \top \ldots \top p_{1})= \sum_{ \sigma \mbox{ (j,l)-shuffle}} p_{\sigma(j+l)} \top \ldots \top p_{\sigma(1)},$$
where $*$ is the product of $gr(\HR)$;
\item[\textnormal{iii)}] $\forall p_j \top \ldots \top p_1 \in Im(F_j)$, 
\begin{eqnarray*}
\Delta(p_j \top \ldots \top p_1 )& =& (1 \otimes p_j \top \ldots \top p_1) + (p_j \top \ldots \top p_1 \otimes 1) \\
&&+ \sum_{k=2}^{k=j} (p_j \top \ldots \top p_{k} ) \otimes (p_{k-1} \top \ldots \top p_1);
\end{eqnarray*}
\item[\textnormal{iv)}] $\forall x \in Im(F_j), j \geq 1, \varepsilon(x)=0;$
\item[\textnormal{v)}] $\forall p_1 \top \ldots \top p_j \in Im(F_j),$ $S_*(p_j \top \ldots \top p_1) = (-1)^j \, p_1 \top \ldots \top p_j.$ 
\end{description}

Clearly, the linear map from $gr(\HR)$ into $\HR$ wich is the identity on every $Im(F_i)$ is a coalgebra isomorphism. It is not an algebra morphism, although we shall prove later that $gr(\HR)$ and $\HR$ are in fact isomorphic Hopf algebras, via another map.

We are going to classify the coalgebra endomorphisms $\HR$ or indifferently $gr(\HR)$.\\
First we fix a notation. Let $u$ be a linear map from $Prim(\HR)^{\otimes i}$ into $Prim(\HR)^{\otimes j}$. Then $\overline{u}$ is the linear map
from $Im(F_i)$ into $Im(F_j)$ defined by $\overline{u}= F_j \circ u \circ  F_i^{-1}$.

\newtheorem{cogebreend}[forme]{Theorem}
\begin{cogebreend}
For all $i \in \mathbb{N}^*$, let $u_i: Prim(\HR)^{\otimes i} \longmapsto Prim(\HR)$. Let $\Phi_{(u_i)}$ be the linear map defined by:
\begin{eqnarray*}
 \Phi_{(u_i)}(1)&=&1;\\
\Phi_{(u_i)}(p_n \top \ldots \top p_1)&=& \sum_{k=1}^n \esp \sum_{a_1 + \ldots +a_k = n} (\overline{ u_{a_1} \otimes \ldots \otimes u_{a_k}}\,)(p_n \top \ldots \top p_1).
\end{eqnarray*}
Then $\Phi_{(u_i)}$ is a coalgebra endomorphism of $\HR$ (or $gr(\HR)$).\\
Moreover, if $\Phi$ is a coalgebra endomorphism of $\HR$ (or $gr(\HR)$), then for all $i \in \mathbb{N}^*$, there exists a unique $ u_i: Prim(\HR)^{\otimes i} \longmapsto Prim(\HR)$,
such that $\Phi = \Phi_{(u_i)}$.
\end{cogebreend} 
\noindent  {\it Proof:} 
first we prove that $\Phi_{(u_i)}$ is a coalgebra endomorphism:
$$\Phi_{(u_i)} \otimes \Phi_{(u_i)} (\tilde{\Delta}(p_n \top \ldots \top p_1))=$$
$$\hspace*{-0.5cm} \sum_{j} \sum_{a_1+ \ldots +a_k  =  j} \mbox{ }
\sum_{b_1 + \ldots +b_l  = n-j} \left[ (\overline{u_{a_1} \otimes \ldots \otimes u_{a_k}}\,) \otimes (\overline{u_{b_1} \otimes \ldots \otimes u_{b_l}}\,) \right]
                \left[ (p_n \top \ldots \top p_{j+1} ) \otimes ( p_{j} \top \ldots \top p_1 ) \right] $$
$$= \tilde{\Delta}\left( \sum_{d_1 + \ldots + d_m = n} (\overline{ u_{d_1} \otimes \ldots \otimes u_{d_m} }\,) ( p_n \top \ldots \top p_1) 
          - \overline{u_n}(p_n \top \ldots \top p_1) \right)$$
$$= \tilde{\Delta}\left( \sum_{d_1 + \ldots + d_m = n} (\overline{ u_{d_1} \otimes \ldots \otimes u_{d_m} }\,) ( p_n \top \ldots \top p_1) 
         \right)-0$$
$$= \tilde{\Delta}(\Phi_{(u_i)}(p_n \top \ldots \top p_1)).$$
Let $\Phi$ be a coalgeabra endomorphism.
$\Delta(\Phi(1)) = \Phi(1) \otimes \Phi(1)$, so $\Phi(1)=0$ or $1$. As $\varepsilon\circ \Phi=\varepsilon$, $\Phi(1)=1$.
We constuct $u_i$ by induction on $i$. For $i=1$, $u_1$ is the restriction of $\Phi$ on $Prim(\HR)$. Suppose we have $u_i$ for $i<n$.
Then with $u'_i=u_i$ if $i<n$ and $u'_i=0$ if $i\geq n$, $\Phi = \Phi_{(u'_i)}$ on $(1)\oplus \oplus_{1}^{n-1}Im(F_j)$.
So 
\begin{eqnarray*}
\Tilde{\Delta}(\Phi(p_n \top \ldots \top p_1))&=& (\Phi \otimes \Phi)\circ \tilde{\Delta}(p_n \top \ldots \top p_1)\\
& =& (\Phi_{(u'_i)} \otimes \Phi_{(u'_i)})\circ \tilde{\Delta}(p_n \top \ldots \top p_1)= \Tilde{\Delta}(\Phi_{(u'_i)}(p_n \top \ldots \top p_1)).
\end{eqnarray*}
So we can take $\overline{u_n}(p_n \top \ldots \top p_1)= (\Phi-\Phi_{(u'_i)})( p_n \top \ldots \top p_1)$.\\
For the uniqueness, observe that $\pi_1 \circ \Phi = \overline{u_i}$ on $Im(F_i)$.\\

We now give a criterion of inversibility of a coalgebra endomorphism:
\newtheorem{inverse}[forme]{Proposition}
\begin{inverse}
\label{invers}
$\Phi_{(u_i)}$ is bijective if and only if the restriction $u_1$ of $\Phi_{(u_i)}$ to $Prim(\HR)$ is bijective. 
\end{inverse}
\noindent  {\it Proof:} $\Rightarrow$: obvious.\\
$\Leftarrow$: we put $\Phi= \Phi_{(u_i)}$.
Recall that $C_i=(1) \oplus \oplus_{1}^{i} Im(F_j)$. As $\Phi(C_i) \subset C_i$, it is enough to show that $\Phi_{\mid C_i}:C_i \longmapsto C_i$ is inversible $\forall i$.
For $i=1$, it is the hypothesis. Suppose it is true for a certain $i-1$. 
Then  $\Phi(p_i \top \ldots \top p_{1})- (\overline{u_1 \otimes \ldots \otimes u_1}\,)(p_i \top \ldots \top p_{1})$ belongs to $C_{i-1}$, so it belongs to $Im(\Phi)$; 
hence $(\overline{u_1 \otimes \ldots \otimes u_1}\,)(C_i) \subset Im(\Phi)$. As $(\overline{u_1 \otimes \ldots \otimes u_1}\,)$ is surjective (because $u_1$ is surjective), $\Phi_{\mid C_i} $ is surjective.\\
Let $x \in C_i$, $\Phi(x)=0$. $x=x_i + y$, $x_i \in Im(F_i)$, $y \in C_{i-1}$. 
Then $\Phi(x)=0= ( \overline{u_1 \otimes \ldots \otimes u_1}\,)(x_i) +C_{i-1}$, so $(\overline{u_1 \otimes \ldots \otimes u_1}\,)(x_i)= 0$ (because it belongs to $Im(F_i)\cap C_{i-1}$). 
As $u_1$ is injective, $x_i=0$, and $x \in C_{i-1}$. As $\Phi_{\mid C_{i-1}}$ is injective, $x=0$: $\Phi_{\mid C_i}$ is injective.\\

 We now give a criterion to know when a coalgebra endomorphism is in fact a bialgebra endomorphism.
\newtheorem{critere2}[forme]{Proposition}
\begin{critere2}
\label{criter2}
Let $\Phi=\Phi_{(u_i)}$ be a coalgebra endomorphism. Let $\Phi^{(n)}=\Phi_{(u_i^n)}$ be the coalgebra endomorphism with
$ u_i^n=u_i$ if $i \leq n$, $u_i^n= 0$ if $i>n$.
\begin{enumerate}

\item (case of $\HR$) $\Phi$ is a bialgebra endomorphism if and only if for all $x_i \in Im(F_i)$, $x_j \in Im(F_j)$,
   $\overline{u_{i+j}}(x_i * x_j ) = -\Phi^{(i+j-1)}(x_i.x_j) + \Phi^{(i+j-1)}(x_i) . \Phi^{(i+j-1)}(x_j)$.
\item (case of $gr(\HR))$ $\Phi$ is a bialgebra endomorphism if and only if for all $x_i \in Im(F_i)$, $x_j \in Im(F_j)$,
   $\overline{u_{i+j}}(x_i * x_j ) = -\Phi^{(i+j-1)}(x_i*x_j) + \Phi^{(i+j-1)}(x_i) * \Phi^{(i+j-1)}(x_j)$.
    
\end{enumerate}
\end{critere2}

\noindent  {\it Proof:} we study the case of $\HR$. 
Observe that $\Phi(x_i.x_j)= \overline{u_{i+j}} ( x_i * x_j) + \Phi^{(i+j-1)}(x_i.x_j)$ because $x_i.x_j - x_i *x_j$ belongs to $ C_{i+j-1}$.
Moreover $ \Phi = \Phi^{(i+j-1)}$ on $C_{i+j-1}.$ 
It is then obvious.
The proof in the case  of $gr(\HR)$ is analog, even easier.

\section{Automorphisms of $\HR$} 

\esp In the following, we shall identify $gr(\HR)$ and $\HR$ as vector spaces via:\\
 $ Id: Im(F_i) \subset gr(\HR) \longmapsto Im(F_i) \subset \HR$.\\
 Now the vector space $\HR$ has two Hopf algebra structures: $(\HR, ., \Delta, S)$ and $(\HR,*,\Delta,S_*)$. Note that the coproduct is the same in both cases.
Both are graded as Hopf algebras by the weight. We still denote by ${\cal H}_i$ the homogeneous components, which are the same for both structures. $(\HR,*,\Delta,S_*)$ is by construction graded as Hopf algebra by $deg_p$, and the homogeneous components are the $Im(F_i)$'s. 

We denote the augmentation ideal, which is  the same for both structures, by $\M$, and its square in $(\HR,.)$  by $\M^2$.
We put $\M_i=\M \cap {\cal H}_i$ and $M_i^2=\M^2 \cap {\cal H}_i$.  
We have: $$\M = \oplus_i \:\M_i \mbox{ and } \M^2=\oplus_i \: \M^2_i. $$
Obviously, $\sum_{j} {\cal H}_i \cap Im(F_j)= {\cal H}_i \cap \sum_j Im(F_j)= {\cal H}_i$ if $ i \leq 1$. So $\M_i^2+ \sum_j {\cal H}_i \cap Im(F_j)= {\cal H}_i=\M_i$. Hence, we can choose $V_{i,j} \subset {\cal H}_i \cap Im(F_j)$, such that $\M_i = \M_i^2 \oplus \oplus_j V_{i,j}$. 
We put $V_i=\oplus_j V_{i,j}$, and $V=\oplus_{i,j} V_{i,j}$.
Note that $V_1={\cal H}_1$. Moreover, for any $x \in \M^2$, $\pi_c(x)=0$, so by lemma \ref{lem4}, $M^2 \cap Im(F_1)= M^2 \cap Prim(\HR)=(0)$. So  $V_{i,1}={\cal H}_i \cap Im(F_1)$.

\newtheorem{sousalgebre}[forme]{Lemma}
\begin{sousalgebre}
\label{sousalg}
$V$ generates the algebra $(\HR,.)$.
\end{sousalgebre}
\noindent {\it Proof:} 
we denote by $\langle V \rangle$ the subalgebra of $(\HR,.)$ generated by $V$.\\
We have to show that ${\cal H}_i \subset \langle V \rangle $ $\forall i \geq 1$. We proceed by induction on $i$. If $i=1$, then it is true since $V_1={\cal H}_1$.
Suppose it is true for any $i' \leq i-1$. Let $x \in {\cal H}_i=\M_i^2 \oplus V_i$.
It is obvious if $x \in V_i$. If $x \in M_i^2$, one can suppose that $x=m_1m_2$, with $m_1$ and $m_2$ in $\M$. Then $m_1$ and $m_2$ cannot be constant, so $weight(m_1)<i$ and $weight(m_2)<i$. So they are in $\langle V \rangle$, so $x \in \langle V \rangle$.

\newtheorem{generateurs}[forme]{Lemma}
\begin{generateurs}
\label{generateur}
$V$ generates the algebra $(\HR,*)$.
\end{generateurs}

\noindent {\it Proof:} 
we denote by $\langle V \rangle_*$ the subalgebra of $(\HR,*)$ generated by $V$.\\
Let $x \in \HR$. Let $j=deg_p(x)$. If $j=1$, then $x \in \langle V \rangle_*$ since $Im(F_1)=\oplus_i V_{i,1}$.
Suppose that $y \in \langle V \rangle_*$ for any $y$ with $deg_p(y)<j$. One can suppose that $x \in \M=\M^2 \oplus V$. If $x \in V$, then $x \in   \langle V \rangle_*$. If $x \in \M^2$, one can suppose that $x=m_1m_2$, with $m_1,m_2 \in \M$.
Then $deg_p(x)=deg_p(m_1)+deg_p(m_2)$, so $deg_p(m_1)<j$ and $deg_p(m_2)<j$, so $m_1$ and $m_2$ are in $\langle V \rangle_*$, and $m_1 *m_2 \in \langle V \rangle_*$. 
By construction of the product $*$, $m_1m_2=m_1*m_2+(1) \oplus \oplus_{k <j} Im(F_k)$. So by induction hypothesis, $x=m_1m_2 \in \langle V \rangle_*$.\\

We denote by $S(V)$ the symmetric algebra generated by $V$.

\newtheorem{dimensionsgene}[forme]{Lemma}
\begin{dimensionsgene}
\label{dimgene}
\begin{enumerate}
\item $\forall i \in \mathbb{N}^*$, $dim(V_i)$ is the number of rooted trees of weight $i$.
\item There is an algebra isomorphism between $(\HR,.)$ and $S(V)$ which is the identity on $V$.
\end{enumerate}
\end{dimensionsgene}

\noindent {\it Proof:} 

1. We have $dim(V_i)= dim(\M_i)-dim(\M_i^2)$. A basis of $\M_i$ is formed by forests of weight $i$, whereas a basis of $\M_i^2$ is formed by non connected forests of weight $i$. The first point is then obvious.

2. As $\HR$ is commutative, we have an algebra morphism:
$$ \Lambda: \left\{
\begin{array}{ccc}
 S(V) & \longmapsto &(\HR,.) \\
  x \in  V & \longmapsto & x 
\end{array}
\right.
$$

By lemma \ref{sousalg},  $\Lambda$ is surjective.
$S(V)$ is graded as algebra by putting $V_i$ in degree $i$. By the first point, the homogeneous components of $S(V)$ and $\HR$ (for the weight) have the same (finite) dimensions. 
Moreover, $\Lambda$ is homogeneous of degree zero; as it is surjective, it is injective; so it is an isomorphism.
\\

Using $\Lambda$, we define an algebra isomorphism:
$$ \Xi: \left\{
\begin{array}{ccc}
 (\HR,.) & \longmapsto &(\HR,*) \\
  x \in  V & \longmapsto & x 
\end{array}
\right.
$$

By lemma \ref{generateur}, $\Xi$ is surjective. Moreover, it is homogenous of degree zero for the weight; 
as the homogeneous components have the same finite dimensions in $(\HR,.)$ and in $(\HR,*)$, it is an isomorphism.

As the coproduct is the same for both Hopf algebra structures on $\HR$, and since $\Xi$ fix a system of generators, it is a bialgebra isomorphism.
Moreover, $\Xi \circ S \circ \Xi^{-1}$ is an antipode of $(\HR,*,\Delta)$, so it is equal to $S_*$. Hence, $\Xi$ is a Hopf algebra isomorphism.

 We have $deg_p(\Xi(x)) \leq deg_p(x)$ $\forall  x \in \HR$, since it is true for any $x \in V$. We get:\\[-3mm]
$$\Xi(\{x \in \HR / deg_p(x)\leq j, weight(x)=i\}) \subset \{x \in \HR / deg_p(x)\leq j, weight(x)=i\} \:\forall i,j.$$
\esp As these spaces have the same finite dimension, they are in fact equal.
We deduce: $$\Xi( \{ x\in \HR / deg_p(x)=j\} )=\{ x\in \HR / deg_p(x)=j\} .$$ \esp We have entirely proved:

\newtheorem{isomorph}[forme]{Theorem}
\begin{isomorph}
\label{isom}
$ gr(\HR)$ and $\HR$ are isomorphic Hopf algebras; there is a Hopf algebra isomorphism $\Xi:(\HR,.) \longmapsto (\HR,*)$ such that $weight(\Xi(x))=weight(x)$ and $deg_p(\Xi(x))=deg_p(x)$ for any $x \in \HR$.

\end{isomorph}
\esp We work now in $gr(\HR)$. 
We denote by $\M^{*2}$ the square of the augmentation ideal in this algebra.
Let $u_1$ be a linear application from $Prim(gr(\HR))$ into itself. Can we extend it to a bialgebra endomorphism of $gr(\HR)$? 
With \ref{criter2}, one sees that $\overline{u_2}$ is entirely determined on $\M^{*2} \cap Im(F_2)$, and we can extend it to the whole $Im(F_2)$ as we want.
More generally, $\overline{u_i}$ is determined over $\M^{*2} \cap Im(F_i)$. So in fact, if we fix a complement $C$ of $\M^{*2}$,
a bialgebra endomorphism $\Phi$ is entirely determined by $(\pi_1 \circ \Phi)_{\mid C}:C \longmapsto Prim(gr(\HR))$. 
Moreover, for any application $L:C \longmapsto Prim(gr(\HR))$, there is a unique bialgebra endomorphism $\Phi_L$ such that $(\pi_1 \circ {\Phi_L})_{\mid C}=L$.     
Because of \ref{isom}, we have the same result for $\HR$. In this case, two important choices of $C$ can be done:
\begin{enumerate}
\item if we choose $C$ the subspace generated by the rooted trees: with notations of \ref{blabla},
$$ (\pi_1 \circ \Phi_{(p_t)})_{\mid C}: 
\left\{
\begin{array}{ccl}
\langle \mbox{ rooted trees } \rangle & \longmapsto &  Prim(\HR)\\
 t & \longmapsto & p_t.
\end{array}
\right.
$$

\item if we choose a complement $C$ which contains $Prim(\HR)$, then we see that \\
$ \left\{
\begin{array}{ccc}
 End_{bialgebra}(\HR) & \longmapsto & {\cal L}(Prim(\HR))\\
 \Phi & \longmapsto & \Phi_{\mid Prim(\HR)}
\end{array}
\right.$
is surjective.
\end{enumerate}

Because of \ref{invers}, we have a surjection:
$$ \chi:
\left\{
\begin{array}{ccl}
 Aut_{bialgeabra}(\HR) & \longmapsto & {GL}(Prim(\HR))\\
 \Phi & \longmapsto & \Phi_{\mid Prim(\HR)}
\end{array}
\right.
$$

 We look for a lifting of ${GL}(Prim(\HR))$ into $Aut_{bialgeabra}(\HR)$. It is easier to work in $gr(\HR)$, 
for there is an obvious lifting: if $u \in {GL}(Prim(\HR))$, we take $u_1=u$, $u_i=0$ if $i\geq 2$;
then one proves easily  that $\Phi_u=\Phi_{(u_i)} \in Aut_{bialgebra}(\HR)$, and $ \Phi_u \circ \Phi_v = \Phi_{u \circ v}$.
So, with \ref{biahopf} we have the following result:
\newtheorem{semidirect}[forme]{Theorem}
\begin{semidirect}
$$ Aut_{bialgebra}(\HR)= Aut_{Hopf}(\HR) = Ker(\chi) \rtimes {GL}(Prim(\HR)).$$
\end{semidirect}
  
\section{Appendix}

$$
\begin{array}{|c|c|c|c|c|c|c|c|c|c|c|c|c|c|c|c|c|c|c|c|c|c|c|c|c|c|c|c|c|c|c|c|}
\hline
n&1&2&3&4&5&6&7&8&9&10&11&12&13&14&15\\
\hline
r_n&1&2&4&9&20&48&115&286&719&1842&4766&12486&32973&87811&235381\\
\hline
h_{n,1}&1&1&1&2&3&8&16&41&98&250&631&1646&4285&11338&30135\\
\hline
\end{array}$$
$$
\begin{array}{|c|c|c|c|c|c|c|c|c|c|c|c|c|c|c|c|c|c|c|c|c|c|c|c|c|c|c|c|c|c|c|c|}
\hline
16&17&18&19&20&21&22&23\\
\hline
634847&1721159&4688676&12826228&35221832&97055181&268282855&743724984\\
\hline
80791&217673&590010&1606188&4392219&12055393&33206321&91752211\\
\hline
\end{array}
$$
$$
\begin{array}{|c|c|c|c|c|c|c|c|c|c|c|c|c|c|c|c|c|c|c|c|c|c|c|c|c|c|c|c|c|c|c|c|}
\hline
24&25&26&27&28&29\\
\hline
2067174645&5759636510&16083734329&45007066269&126186554308&354426847597\\
\hline
254261363&706465999&1967743066&5493195530&15367129299&43073007846\\
\hline
\end{array}
$$


\begin{thebibliography}{9}
\bibitem{Kreimer1}{D. Kreimer}: {\it On the Hopf algebra structure of pertubative quantum field theories}, (1998), q-alg/9707029.
\bibitem{Moscovici}{A. Connes, H. Moscovici}: {\it Hopf algebras, cyclic Cohomology and the transverse Index Theorem}, IHES/M/98/37, math.DG/9806109. 
\bibitem{Connes}{A. Connes, D. Kreimer}: {\it Hopf algebras, Renormalization and Noncommutative geometry}, {\it Comm. Math. Phys} \textbf{199} (1998) 203, hep-th/9808042.
\bibitem{Broadhurst}{D. J. Broadhurst, D. Kreimer}: {\it Renormalization automated by Hopf algebra}, (1998), hep-th/9810087.
\bibitem{Kreimer2}{D. Kreimer}: {\it On Overlapping Divergences}, (1999), hep-th/9810022.
\bibitem{Kreimer} {D. J. Broadhurst, D. Kreimer:} {\it Towards cohomology of renormalization: bigrading the combinatorial Hopf algebra of rooted trees}, (2000), hep-th/0001202.
\bibitem{Panaite}{Panaite:} {\it Relating the Connes-Kreimer and Grossman-Larsen Hopf algebras built on  rooted trees}, (2000), math.QA/0003074.
\bibitem{Kreimer3}{D. Kreimer:} {\it Combinatorics of (pertubative) Quantum Field Theory}, (2000), hep-th/0010059.
\bibitem{Sloane}{N. J. A. Sloane:} {\it On-line Encyclopedia of Integer Sequences}, sequence A000081, \textbf{http://www.research.att.com/$\tilde{\hspace{2mm}}$njas/sequences}. 
\bibitem{Bourbaki}{N. Bourbaki}: {\it Groupes et alg\`ebres de Lie, chapitres 2 et 3}, {(1972)}, {Hermann}, {ch II, $\S 3$, corollaire 2}.
\bibitem{Sweedler}{M. Sweedler}: {\it Hopf algebras}, (1969), {W. A. Benjamin, Inc., New York.}  
\end{thebibliography}
\end{document}